\documentclass[preprint,12pt, sort&compress,numbers]{elsarticle}

\usepackage{soul}
\input{packages.tex}

\renewcommand{\Re}{\mathrm{Re}}

\newcommand{\cO}{\mathcal{O}}

\newcommand{\cT}{\mathcal{T}}

\newcommand{\cQ}{\mathcal{Q}}

\newcommand{\rR}{\mathbb{R}}
\newcommand{\rI}{\mathbb{I}}

\newcommand{\rU}{\mathbb{U}}
\newcommand{\rP}{\mathbb{P}}
\newcommand{\rC}{\mathbb{C}}

\newcommand{\bB}{\mathbf{B}}
\newcommand{\bC}{\mathbf{C}}

\newcommand{\bL}{\mathbf{L}}
\newcommand{\bS}{\mathbf{S}}
\newcommand{\bV}{\mathbf{V}}
\newcommand{\bU}{\mathbf{U}}
\newcommand{\bP}{\mathbf{P}}

\newcommand{\bK}{\mathbf{K}}
\newcommand{\bI}{\mathbf{I}}
\newcommand{\bZ}{\mathbf{Z}}

\newcommand{\bn}{\mathbf{n}}
\newcommand{\bq}{\mathbf{q}}
\newcommand{\bx}{\mathbf{x}}
\newcommand{\bu}{\mathbf{u}}

\newcommand{\bw}{\mathbf{w}}
\newcommand{\bus}{\bu^{\dagger}}
\newcommand{\bps}{\bp^{\dagger}}
\newcommand{\bg}{\mathbf{g}}

\newcommand{\tbu}{\tilde{\bu}}
\newcommand{\tp}{\tilde{p}}
\newcommand{\bqs}{\mathbf{q}^{\dagger}}

\newcommand{\bp}{\mathbf{p}}
\newcommand{\f}{\mathbf{f}}
\newcommand{\bPhi}{\mathbf{\Phi}}
\newcommand{\bphi}{\bm{\phi}}

\newcommand{\ps}{p^{\dagger}}

\newcommand{\hL}{{\hat{L}}}
\newcommand{\hB}{{\hat{B}}}
\newcommand{\hK}{{\hat{K}}}
\newcommand{\hC}{{\hat{C}}}

\newcommand{\hu}{{\hat{u}}}
\newcommand{\hus}{{\hat{u}^{\dagger}}}
\newcommand{\hbq}{{\hat{\bq}}}

\newcommand{\hbu}{{\hat{\bu}}}
\newcommand{\hbus}{{\hat{\bu}^{\dagger}}}
\newcommand{\hbp}{{\hat{\bp}}}
\newcommand{\hbps}{{\hat{\bp}^{\dagger}}}
\newcommand{\hbB}{{\hat{\bB}}}
\newcommand{\hbC}{{\hat{\bC}}}

\newcommand{\hbK}{{\hat{\bK}}}
\newcommand{\hbL}{{\hat{\bL}}}

\newcommand{\uhbC}{\underline{\hbC}}

\newcommand{\bOmega}{\overline{\Omega}}

\def\inprod#1#2{\left\langle #1, #2 \right\rangle}
\def\form#1#2#3{#1\!\inprod{#2}{#3}}
\def\form#1#2#3{#1\!\inprod{#2}{#3}}
\def\form#1#2#3{#1\!\inprod{#2}{#3}}

\def\jump#1{\left\llbracket #1 \right\rrbracket}
\def\avg#1{\left\{\!\!\left\{ #1 \right\}\!\!\right\}}

\newcommand{\trun}{\mathrm{Tr}}

\newcommand{\logLogSlopeTriangle}[6]
{

    \pgfplotsextra
    {
        \pgfkeysgetvalue{/pgfplots/xmin}{\xmin}
        \pgfkeysgetvalue{/pgfplots/xmax}{\xmax}
        \pgfkeysgetvalue{/pgfplots/ymin}{\ymin}
        \pgfkeysgetvalue{/pgfplots/ymax}{\ymax}

        \pgfmathsetmacro{\xArel}{#1}
        \pgfmathsetmacro{\yArel}{#3}
        \pgfmathsetmacro{\xBrel}{#1-#2}
        \pgfmathsetmacro{\yBrel}{\yArel}
        \pgfmathsetmacro{\xCrel}{\xArel}

        \pgfmathsetmacro{\lnxB}{\xmin*(1-(#1-#2))+\xmax*(#1-#2)} 
        \pgfmathsetmacro{\lnxA}{\xmin*(1-#1)+\xmax*#1} 
        \pgfmathsetmacro{\lnyA}{\ymin*(1-#3)+\ymax*#3} 
        \pgfmathsetmacro{\lnyC}{\lnyA+#4*(\lnxA-\lnxB)}
        \pgfmathsetmacro{\yCrel}{(\lnyC-\ymin)/(\ymax-\ymin)} 

        \coordinate (A) at (rel axis cs:\xArel,\yArel);
        \coordinate (B) at (rel axis cs:\xBrel,\yBrel);
        \coordinate (C) at (rel axis cs:\xCrel,\yCrel);

        \draw[#5]   (A)-- node[pos=0.5,anchor=#6] {1}
                    (B)-- 
                    (C)-- node[pos=0.5,anchor=west] {#4}
                    cycle;
    }
}

\newcommand{\logLogReverseSlopeTriangle}[6]
{

    \pgfplotsextra
    {
        \pgfkeysgetvalue{/pgfplots/xmin}{\xmin}
        \pgfkeysgetvalue{/pgfplots/xmax}{\xmax}
        \pgfkeysgetvalue{/pgfplots/ymin}{\ymin}
        \pgfkeysgetvalue{/pgfplots/ymax}{\ymax}

        \pgfmathsetmacro{\xArel}{#1}
        \pgfmathsetmacro{\yArel}{#3}
        \pgfmathsetmacro{\xBrel}{#1+#2}
        \pgfmathsetmacro{\yBrel}{\yArel}
        \pgfmathsetmacro{\xCrel}{\xArel}

        \pgfmathsetmacro{\lnxB}{\xmin*(1-(#1-#2))+\xmax*(#1-#2)} 
        \pgfmathsetmacro{\lnxA}{\xmin*(1-#1)+\xmax*#1} 
        \pgfmathsetmacro{\lnyA}{\ymin*(1-#3)+\ymax*#3} 
        \pgfmathsetmacro{\lnyC}{\lnyA-#4*(\lnxA-\lnxB)}
        \pgfmathsetmacro{\yCrel}{(\lnyC-\ymin)/(\ymax-\ymin)} 

        \coordinate (A) at (rel axis cs:\xArel,\yArel);
        \coordinate (B) at (rel axis cs:\xBrel,\yBrel);
        \coordinate (C) at (rel axis cs:\xCrel,\yCrel);

        \draw[#5]   (A)-- node[pos=0.5,anchor=#6] {1}
                    (B)-- 
                    (C)-- node[pos=0.5,anchor=east] {#4}
                    cycle;
    }
}

\pgfplotsset{select coords between index/.style 2 args={
    x filter/.code={
        \ifnum\coordindex<#1\fi
        \ifnum\coordindex>#2\fi
    }
}}

\newif\ifcommentcoloring
\commentcoloringfalse
\newif\ifrefereecoloring
\refereecoloringfalse

\newcommand{\YC}[1]{\ifcommentcoloring \textcolor{red}{#1}\else #1\fi}
\newcommand{\PR}[1]{\ifcommentcoloring \textcolor{blue}{#1}\else #1\fi}
\newcommand{\KC}[1]{\ifcommentcoloring \textcolor{brown!70!black}{#1}\else #1\fi}
\newcommand{\TYL}[1]{\ifcommentcoloring \textcolor{magenta}{#1}\else #1\fi}

\journal{Journal of Computational Physics}

\begin{document}

\begin{frontmatter}



\title{
Scaled-up prediction of steady Navier-Stokes equation with component reduced order modeling
}


\author[casc]{Seung Whan Chung\corref{cor1}}
\cortext[cor1]{corresponding author}
\ead{chung28@llnl.gov}
\author[casc]{Youngsoo Choi}
\author[aeed]{Pratanu Roy}
\author[ced]{Thomas Roy}
\author[ced]{Tiras Y. Lin}
\author[med]{Du T. Nguyen}
\author[leaf]{Christopher Hahn}
\author[cemm]{Eric B. Duoss}
\author[msd]{Sarah E. Baker}

\affiliation[casc]{organization={Center for Applied Scientific Computing, Lawrence Livermore National Laboratory},
  city={Livermore},
  postcode={94550}, 
  state={CA},
  country={US}}
\affiliation[aeed]{organization={Atmospheric, Earth and Energy Division, Lawrence Livermore National Laboratory},
  city={Livermore},
  postcode={94550}, 
  state={CA},
  country={US}}
\affiliation[ced]{organization={Computational Engineering Division, Lawrence Livermore National Laboratory},
  city={Livermore},
  postcode={94550}, 
  state={CA},
  country={US}}
\affiliation[msd]{organization={Material Science Division, Lawrence Livermore National Laboratory},
  city={Livermore},
  postcode={94550}, 
  state={CA},
  country={US}}
\affiliation[leaf]{organization={Laboratory for Energy Applications for the Future, Lawrence Livermore National Laboratory},
  city={Livermore},
  postcode={94550}, 
  state={CA},
  country={US}}
\affiliation[cemm]{organization={Center for Engineered Materials, Manufacturing and Optimization, Lawrence Livermore National Laboratory},
  city={Livermore},
  postcode={94550}, 
  state={CA},
  country={US}}

\begin{abstract}
Scaling up new scientific technologies from laboratory to industry often involves demonstrating performance on a larger scale.
Computer simulations can accelerate design and predictions in the deployment process,
though traditional numerical methods are computationally intractable even for intermediate pilot plant scales.
Recently, component reduced order modeling method is developed to tackle this challenge
by combining projection reduced order modeling and discontinuous Galerkin domain decomposition.
However, while many scientific or engineering applications involve nonlinear physics,
this method has been only demonstrated for various linear systems.
In this work, the component reduced order modeling method is extended
to steady Navier-Stokes flow, with application to general nonlinear physics in view.
Large-scale, global domain is decomposed into combination of small-scale unit component.
Linear subspaces for flow velocity and pressure are identified via proper orthogonal
decomposition over sample snapshots collected at small scale unit component.
Velocity bases are augmented with pressure supremizer, in order to satisfy inf-sup condition for stable pressure prediction.
Two different nonlinear reduced order modeling methods are employed and compared for efficient evaluation of nonlinear advection: 3rd-order tensor projection operator and empirical quadrature procedure.
The proposed method is demonstrated on flow over arrays of five different unit objects,
achieving $23$ times faster prediction with less than $4\%$ relative error up to $256$ times larger scale domain than unit components.
Furthermore, a numerical experiment with pressure supremizer strongly indicates the need of supremizer for stable pressure prediction. A comparison between tensorial approach and empirical quadrature procedure is performed, which suggests a slight advantage for empirical quadrature procedure.
\end{abstract}


\begin{highlights}
\item A novel component reduced order model with bottom-up training is developed for scaled-up, robust and accelerated prediction of steady Navier-Stokes flow.
\item Tensorial method and empirical quadrature procedure are employed and compared for
reduced order modeling of nonlinear advection.
\item Bases for velocity linear subspace are augmented with pressure supremizer
and its impact on stable pressure prediction is demonstrated.
\item The constructed reduced order model can provide sufficiently accurate, robust predictions at large scales in an efficient manner.
\end{highlights}


\end{frontmatter}


\section{Introduction}




\PR{Scaling up new scientific technologies from laboratory to industry often requires validating their effectiveness at an industrial scale.}
Traditionally, this is done with physical pilot plants, which are expensive and time-consuming.
\PR{Although simulations using traditional numerical methods can accelerate design and predictions, conducting simulations at pilot scales, that are much larger than lab scales, is usually resource intensive and may be impractical. To obtain fast predictions at the pilot scale, simplified volume-averaged models are frequently employed. However, this approach sacrifices crucial physical phenomena that take place at smaller scales.}

\par
The recently proposed component reduced order modeling (CROM)~\cite{chung2024train} tackles this challenge by combining projection-based reduced order modeling (PROM)
with discontinuous Galerkin domain decomposition (DG-DD).
From small scale physics solution data,
Proper orthogonal decomposition (POD)~\cite{berkooz1993proper} can extract a low-dimensional linear subspace
that can effectively represent the physics solutions at small scales.
\PR{This} small-scale reduced-order model (ROM) solves the physics governing equation on the reduced linear subspace,
thereby achieving both \PR{desired} accuracy and cheap computation time.
Multiple small-scale unit ROMs can then be assembled into a large-scale global ROM,
where the interface condition is handled via discontinuous Galerkin penalty terms.
Such combination of PROM and DG-DD enables robust scaled-up predictions solely with the ROMs at small unit scales.
Furthermore, CROM provides great flexibility and simplicity as the interface condition is handled at the full-order model (FOM) level
and no additional handling is required at the ROM level.
\par
CROM has been successfully demonstrated for several linear physics problems such as
the Poisson equation, Stokes flow, Helmholtz partial differential equations, linear elasticity, and the advection-diffusion equation~\cite{chung2024train,huynh2013static,eftang2014port,vallaghe2014static,eftang2013port,smetana2016optimal,mcbane2021component,mcbane2022stress}.
However,
\TYL{in addition to the added computational cost required for scaled up systems, the equations needed to accurately capture the physical phenomena likely are nonlinear. These nonlinearities could arise directly from the form of the governing equation derived from first principles, or from correlations fit to experimental data that are used to effectively describe some complex phenomena. For example, in carbon capture columns, where CO$_2$-laden gas flows against a CO$_2$ absorbing solvent through a packed column, engineering correlations with complex mathematical forms \citep{onda1968mass,lin2024advancing} describe the wetted surface area of the packing material. \PR{The liquid solvent distribution over the packing column depends on the complex interactions of gas-liquid interfacial velocities, surface tension, and contact angles~\citep{Singh2022,ellebracht20233d}}. 
}
\KC{
In CO$_2$ electrolysis, reaction rates are often described as a power law combined with an exponential term \citep{ehlinger2024modeling}.
}
\TYL{In homogenized porous media flow, pressure drops are  linear in the flow rate when the flow is weak (the Darcy equation), but can become quadratic in the flow rate as the flow rate rate is increased (the Ergun equation). Additionally, in flowing systems, increases in scale are typically related to corresponding increases in the Reynolds number, thus requiring the inclusion of the nonlinear advection term. Thus, any tool that can aid in the scale up of simulations should also be able to handle nonlinear problems to be most impactful.}
\par
\YC{Several CROMs have been extended to address nonlinear problems. Hoang et al. \cite{hoang2021domain} applied the least-squares Petrov–Galerkin (LSPG) formulation along with hyper-reduction techniques to efficiently manage nonlinear terms within the CROM framework. However, their method requires separately tracking the port basis from the interior snaptshot data, which can be cumbersome from an implementation standpoint. More recently, Diaz et al. \cite{diaz2024fast} advanced Hoang's CROM work by incorporating a nonlinear manifold solution representation, though their training process followed a top-down approach rather than a bottom-up one. Ebrahimi and Yano \cite{ebrahimi2024hyperreduced} introduced a hyper-reduced CROM that enables bottom-up training and addresses nonlinear problems. Nevertheless, their approach still necessitates the tracking of an interface basis, commonly referred to as the port basis.}
\KC{
Wentland~\textit{et al.}~\cite{wentland2024role} have developed a domain-decomposed nonlinear ROM equipped with Schwarz alternating method.
Due to the nature of Schwarz alternating method, the domain should be decomposed with numerically overlapping interface.
Furthermore, their training process follows a top-down approach and does not envision a bottom-up one.
}
\par
In this work, we extend CROM to the steady incompressible Navier-Stokes equations \YC{without the need for a separate interface basis,}
with potential extension to general nonlinear systems in mind.
In general, a naive projection of the nonlinear term onto a linear subspace would not gain any speed-up,
as it cannot be precomputed and requires evaluation at all grid points.
For the advection term in the Navier-Stokes equations,
the tensorial approach has been widely used~\cite{Lassila2014}, exploiting the fact that the advection is quadratic in terms of velocity.
Hyper-reduction is an alternative nonlinear ROM approach,
where the projection of the nonlinear term is approximated with a sparsely sampled grid evaluation~\cite{Willcox2006,Chaturantabut2010,chapman2017accelerated,Yano2019dg}.
Among various hyper-reduction approaches,
the empirical quadrature procedure (EQP)~\cite{chapman2017accelerated,Yano2019dg} is designed specifically for weak-form approximation
and is favorable for finite-element type discretizations and PROM.
In this work, the tensorial approach and EQP method are employed and compared for their accuracy and computational cost.
\par
Another challenge arises due to the combined effect of incompressibility and nonlinearity of the Navier-Stokes equations.
While, in Stokes flow, pressure can be represented as a linearly dependent variable of velocity~\cite{chung2024train},
such representation is not feasible for Navier-Stokes flow and velocity and pressure must be treated as independent variables.
However, the velocity linear subspace identified from incompressible flow solutions also remains incompressible,
violating the uniqueness condition for pressure in incompressible flow systems, also known as the inf-sup condition~\cite{Babuvska1971,Brezzi1974,Ladyzhenskaya1963,Taylor1973}.
While the FOM finite element space may have already satisfied this condition,
the resulting ROM solution space does not necessarily satisfy the same necessary physics condition,
leading to spurious pressure predictions.
This issue can be addressed by augmenting the velocity linear subspace with compressible velocity components.
Ballarin \textit{et al.}~\cite{ballarin2015supremizer} introduced supremizer enrichment procedure,
which constructs such compressible components from gradients of pressure POD modes.
In this study, we employ the supremizer enrichment in the context of CROM.
In essence, the velocity basis of each unit reference component is augmented with supremizers from its respective pressure POD basis.
\par
The rest of the paper is organized as follows. In Section~\ref{sec:formulation}, we provide a detailed formulation of the proposed
component model reduction approach for the steady incompressible Navier-Stokes equations.
The proposed method is then demonstrated in Section~\ref{sec:results}
for a scaled-up prediction of flow past an array of objects at a moderate Reynolds number.
Lastly, we conclude this paper in Section~\ref{sec:conclusion} with a discussion of future works.
\section{General framework}\label{sec:formulation}

\subsection{Governing physics equations}
In the steady case, the incompressible Navier--Stokes equations read
\begin{subequations}
    \begin{equation}
        -\nu\nabla^2\tbu + \nabla\tp + \tbu\cdot\nabla\tbu = 0
    \end{equation}
    \begin{equation}
        \nabla\cdot\tbu = 0,
    \end{equation}
\end{subequations}
for flow velocity $\tilde{\bu}\in H^1(\Omega)^{d}$, and pressure $\tilde{p} \in H^1(\Omega)$ on a global domain $\Omega \subset \rR^{d}$ with $d=2, 3$.
The non-dimensional viscosity $\nu=1/\mathrm{Re} = \frac{\tilde{\nu}_0}{\tilde{U}_0\tilde{L}_0}$ is the inverse of Reynolds number,
which is defined based on a nominal inflow velocity $\tilde{U}_0$, length scale $\tilde{L}_0$, and kinematic viscosity of fluid $\tilde{\nu}_0$.
\par
The boundary of the global domain $\partial\Omega$ is composed of Dirichlet boundary $\partial\Omega_{di}$ and Neumann boundary $\partial\Omega_{ne}$,
i.e., $\partial\Omega = \partial\Omega_{di}\cup\partial\Omega_{ne}$.
The velocity Dirichlet boundary condition is
\begin{equation}
\tilde{\bu} = \bg_{di} \qquad \text{on } \partial\Omega_{di},
\end{equation}
with a prescribed flow velocity $\bg_{di}$.
For $\partial\Omega_{ne}$,
we consider a homogeneous Neumann condition for outgoing flow,
\begin{equation}
\bn\cdot(-\nu\nabla\tilde{\bu} + \tilde{p}\bI) = \mathbf{0} \qquad \text{on } \partial\Omega_{ne},
\end{equation}
with $\bn$ the outward normal vector of the domain $\Omega$.


\subsection{Domain decomposition}
We consider the global-scale domain $\Omega$ decomposed into $M$ subdomains $\Omega_m$,
\begin{equation}
    \Omega = \bigcup\limits_{m=1}^M\Omega_m,
\end{equation}
where all subdomains can be categorized into a few reference domains,
\begin{equation}\label{eq:ref-domain}
    \Omega_m \in \rC\equiv\left\{\bOmega_1, \bOmega_2, \ldots\right\}
    \qquad \forall \: m=1,\ldots,M.
\end{equation}
The global flow velocity $\tbu$ and $\tp$ are composed of subdomain states $\tbu_m \in H^1(\Omega_m)^{d}$ and $\tp_m \in H^1(\Omega_m)$, respectively,
\begin{equation}
    \tbu = \{\tbu_m\}_{m=1}^M \qquad\text{and}\qquad\tp = \{\tp_m\}_{m=1}^M.
\end{equation}
Each $\tbu_m$ and $\tp_m$ satisfies the governing equations on its subdomain,
\begin{subequations}\label{eq:dd-gov}
    \begin{equation}
        -\nu\nabla^2\tbu_m + \nabla\tp_m + \tbu_m\cdot\nabla\tbu_m = 0
    \end{equation}
    \begin{equation}\label{eq:dd-gov-div}
        \nabla\cdot\tbu_m = 0,
    \end{equation}
    with the continuity and smoothness constraints on the interface $\Gamma_{m,n}\equiv\partial\Omega_m\cap\partial\Omega_n$,
    \begin{equation}
        \jump{\tbu} \equiv \tbu_m - \tbu_n = 0 \qquad \text{on } \Gamma_{m,n}
    \end{equation}
    \begin{equation}
        \jump{\tp} \equiv \tp_m - \tp_n = 0 \qquad \text{on } \Gamma_{m,n}
    \end{equation}
    \begin{equation}
        \avg{\bn\cdot\nabla\tbu} \equiv \frac{1}{2}\left( \bn_m\cdot\nabla\tbu_m + \bn_n\cdot\nabla\tbu_n \right) = 0 \qquad \text{on } \Gamma_{m,n}
    \end{equation}
    \begin{equation}
        \avg{\bn\cdot\nabla\tp} \equiv \frac{1}{2}\left( \bn_m\cdot\nabla\tp_m + \bn_n\cdot\nabla\tp_n \right) = 0 \qquad \text{on } \Gamma_{m,n},
    \end{equation}
\end{subequations}
where $\bn_m$ is the outward normal vector of the subdomain $\Omega_m$,
and $\bn_m = -\bn_n$ on $\Gamma_{m,n}$.
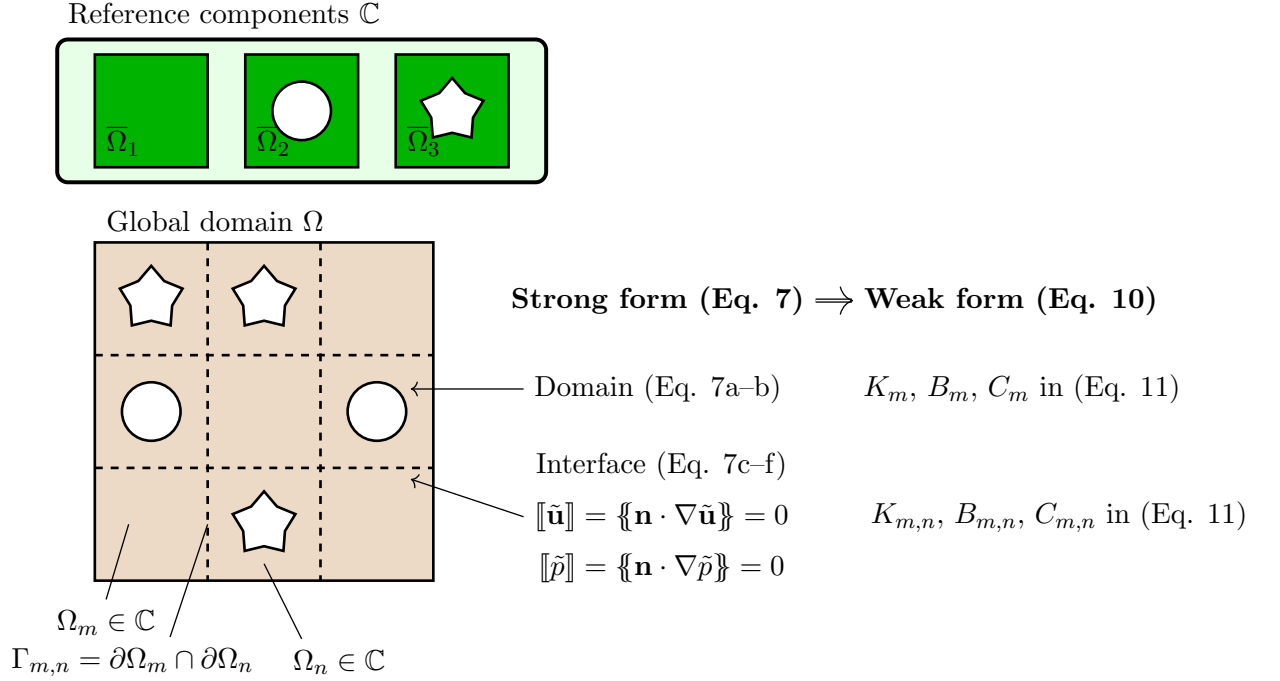
\begin{figure}[tbh]
    \begin{tikzpicture}[font=\small,]

\def\compsize{1.5}
\def\compgap{2.0}

\draw[rounded corners, draw=black, fill=green!10, line width=1.5] (-0.5, -0.2) rectangle (6.0, 1.7);
\node[anchor=south west,] at (-0.5, 1.7) {Reference components $\rC$};

\draw[draw=black, fill=green!70!black, line width=1.0,] (0, 0) node[anchor=south west,] {$\bOmega_1$} rectangle ++ (\compsize, \compsize);

\draw[draw=black, fill=green!70!black, line width=1.0,] (\compgap, 0) node[anchor=south west,] {$\bOmega_2$} rectangle ++ (\compsize, \compsize);
\node[draw=black, fill=white, line width=1.0, circle, scale=2.0,] at (\compgap+0.5*\compsize, 0.5*\compsize) {}; 

\draw[draw=black, fill=green!70!black, line width=1.0,] (2*\compgap, 0) node[anchor=south west,] {$\bOmega_3$} rectangle ++ (\compsize, \compsize);
\node[draw=black, fill=white, line width=1.0, star, scale=1.5,] at (2*\compgap+0.5*\compsize, 0.5*\compsize) {}; 

\draw[draw=black, fill=brown!30, line width=1.0,] (0, -1)
    node[anchor=south west,] {Global domain $\Omega$} rectangle ++(3*\compsize,-3*\compsize);
\draw[draw=black, dashed, line width=1.0,] (\compsize, -1) -- ++(0, -3*\compsize);
\draw[draw=black, dashed, line width=1.0,] (2*\compsize, -1) -- ++(0, -3*\compsize);
\draw[draw=black, dashed, line width=1.0,] (0, -1-\compsize) -- ++(3*\compsize, 0);
\draw[draw=black, dashed, line width=1.0,] (0, -1-2*\compsize) -- ++(3*\compsize, 0);

\node[draw=black, fill=white, line width=1.0, star, scale=1.5,] at (0 + 1.5*\compsize, -1 - 0.5*\compsize) {};
\node[draw=black, fill=white, line width=1.0, star, scale=1.5,] at (0 + 1.5*\compsize, -1 - 2.5*\compsize) {};
\node[draw=black, fill=white, line width=1.0, star, scale=1.5,] at (0 + 0.5*\compsize, -1 - 0.5*\compsize) {};
\node[draw=black, fill=white, line width=1.0, circle, scale=2.,] at (0 + 2.5*\compsize, -1 - 1.5*\compsize) {};
\node[draw=black, fill=white, line width=1.0, circle, scale=2.,] at (0 + 0.5*\compsize, -1 - 1.5*\compsize) {};

\draw[draw=black, line width=0.5,] (0 + 1.5*\compsize, -1 - 2.85*\compsize) -- ++(0.5, -1.0)
    node[anchor=north, xshift=.5cm,] {$\Omega_n\in\rC$};
\draw[draw=black, line width=0.5,] (0 + 0.3*\compsize, -1 - 2.5*\compsize) -- ++(-0.3, -1.0)
    node[anchor=north, xshift=.0cm,] {$\Omega_m\in\rC$};
\draw[draw=black, line width=0.5,] (0 + 1*\compsize, -1 - 2.5*\compsize) -- ++(-0.5, -1.5)
    node[anchor=north, xshift=-.5cm,] {$\Gamma_{m,n}=\partial\Omega_m\cap\partial\Omega_n$};

\draw[<-, line width=0.5,] (0 + 2.8*\compsize, -1 - 1.3*\compsize) -- ++(1.5, 0)
    node[anchor=west,] (domstreq) {Domain (Eq. \ref{eq:dd-gov}a--b)};
\node[anchor=south, yshift=0.5cm,] (streq) at (domstreq.north) {\bf Strong form (Eq. \ref{eq:dd-gov})};
\draw[<-, line width=0.5,] (0 + 2.8*\compsize, -1 - 2.1*\compsize) -- ++(1.5, -0.5)
    node[anchor=west,] (ifstreq) {$\jump{\tbu} = \avg{\bn\cdot\nabla\tbu} = 0$};
\node[anchor=south,] at (ifstreq.north) {Interface (Eq. \ref{eq:dd-gov}c--f)};
\node[anchor=north,] at (ifstreq.south) {$\jump{\tp} = \avg{\bn\cdot\nabla\tp} = 0$};

\draw[->, double, line width=0.5,] (streq.east) -- ++(0.5, 0) node[anchor=west,] {\bf Weak form (Eq. \ref{eq:weak-fom-gov})};
\node[anchor=west, xshift=0.8cm,] at (domstreq.east) {$K_m$, $B_m$, $C_m$ in (Eq. \ref{eq:fom-global-op})};
\node[anchor=west, xshift=0.8cm,] at (ifstreq.east) {$K_{m,n}$, $B_{m,n}$, $C_{m,n}$ in (Eq. \ref{eq:fom-global-op})};

\end{tikzpicture}
    \caption{\PR{Domain decomposition showing the global domain, subdomains and reference components.}}
    \label{fig:dd-illustration}
\end{figure}
Figure~\ref{fig:dd-illustration} illustrates (\ref{eq:dd-gov})
with an example of 3-by-3 global domain constructed with 3 types of components.

\subsection{Full order model with Discontinuous Galerkin (DG)}\label{subsec:dg-dd}

Several choices of finite element spaces for incompressible flow are available
for our domain decomposition framework~\cite{Toselli2002,Cockburn2002,Elman2014}.
Our choice of the finite element space $\rU_{s+1}\otimes\rP_s$ is Taylor-Hood elements~\cite{Taylor1973}.
$\rU_s \subset H^1(\Omega)^d$ and $\rP_s \subset H^1(\Omega)$ are the piecewise $s$-th order polynomial function spaces
composed of $M$ subdomain finite element spaces,
\begin{equation}\label{eq:fes}
    \rU_s = \prod_{m=1}^M \rU_{m, s}\qquad\text{and}\qquad\rP_s = \prod_{m=1}^M \rP_{m, s}.
\end{equation}
The subdomain finite element spaces for velocity and pressure are defined as follows, respectively,
\begin{subequations}
    \begin{equation}
        \rU_{m, s} = \left\{ \bu_m \in H^1(\Omega_m)^d \;\bigg|\; \bu_m\big|_{\kappa} \in V_s(\kappa)^{d} \quad \forall \kappa\in\cT(\Omega_m) \right\}
    \end{equation}
    \begin{equation}
        \rP_{m, s} = \left\{ p_m \in H^1(\Omega_m) \cap H^1_0(\Omega) \;\bigg|\; p_m\big|_{\kappa} \in V_s(\kappa) \quad \forall \kappa\in\cT(\Omega_m) \right\},
    \end{equation}
\end{subequations}
where $\cT(\Omega_m)$ is the set of mesh elements in subdomain $\Omega_m$,
and $V_s(\kappa)$ is the space of $s$th order polynomials in a mesh element $\kappa$.
For Taylor-Hood elements, the velocity space is defined one-order higher than the pressure space
in order to satisfy the Ladyzhenskaya-Babu\v{s}ka-Brezzi (LBB) condition for uniqueness of the solution~\cite{Babuvska1971,Brezzi1974,Ladyzhenskaya1963,Taylor1973}. 
In the case of $\partial\Omega_{ne}=\varnothing$, the pressure is restricted to have mean value zero in the global domain $\Omega$.
\par
Standard Galerkin finite-element discretization seeks approximate solutions $\bu\in\rU_{s+1}$ and $p\in\rP_s$
that satisfy (\ref{eq:dd-gov}) in a weak sense,
\begin{subequations}\label{eq:weak-fom-gov}
    \begin{equation}
        \form{K}{\bus}{\bu} + \form{B}{p}{\bus} + \form{C}{\bus}{\bu} = \form{L}{\bus}{\f} + \form{L_u}{\bus}{\bg}
    \end{equation}
    \begin{equation}
        \form{B}{\ps}{\bu} = \form{L_p}{\ps}{\bg},
    \end{equation}
\end{subequations}
for all test function $\bus\in\rU_{s+1}$ and $\ps\in\rP_s$, where
the functionals $K$, $B$ and $C$ in (\ref{eq:weak-fom-gov}) approximate the viscous flux, velocity divergence and nonlinear advection term, respectively.
They are composed of subdomain, interface and boundary functionals,
\begin{subequations}\label{eq:fom-global-op}
    \begin{equation}
        \form{K}{\bus}{\bu} =
        \sum\limits_{m=1}^{M}\form{K_m}{\bus}{\bu} + 
        \sum\limits_{(m,n)\in \rI}\form{K_{m,n}}{\bus}{\bu}
        + \sum\limits_{m=1}^{M}\form{K_{m,di}}{\bus}{\bu}
    \end{equation}
    \begin{equation}
        \form{B}{\ps}{\bu} =
        \sum\limits_{m=1}^{M}\form{B_m}{\ps}{\bu} + 
        \sum\limits_{(m,n)\in \rI}\form{B_{m,n}}{\ps}{\bu}
        + \sum\limits_{m=1}^{M}\form{B_{m,di}}{\ps}{\bu}
    \end{equation}
    \begin{equation}
        \form{C}{\bus}{\bu} =
        \sum\limits_{m=1}^{M}\form{C_m}{\bus}{\bu} + 
        \sum\limits_{(m,n)\in \rI}\form{C_{m,n}}{\bus}{\bu}
        + \sum\limits_{m=1}^{M}\form{C_{m,di}}{\bus}{\bu}.
    \end{equation}
\end{subequations}
The subdomain functionals $K_m$, $B_m$ and $C_m$ are defined on $\Omega_m$,
the interface functionals $K_{m,n}$, $B_{m,n}$ and $C_{m,n}$ on all interfaces $\Gamma_{m,n}$ with $\rI=\left\{(m,n) \,\middle|\, \Gamma_{m,n}\ne\varnothing\right\}$,
and the Dirichlet boundary functionals $K_{m,di}$, $B_{m,di}$ and $C_{m, di}$ on $\Gamma_{m,di} \equiv \partial\Omega_m \cap \partial\Omega_{di}$.
\par
For linear operators of the viscous flux and velocity divergence,
we follow the formulation of Toselli~\cite{Toselli2002}.
The bilinear form $K$ approximates the viscous term,
\begin{subequations}\label{eq:Kform}
    \begin{equation}
        \form{K_m}{\bus}{\bu} = \inprod{\nu\nabla\bus_m}{\nabla\bu_m}_{\Omega_m}
    \end{equation}
    \begin{equation}
        \begin{split}
            \form{K_{m,n}}{\bus}{\bu} =&
            - \inprod{\avg{\bn\cdot\nu\nabla\bus}}{\jump{\bu}}_{\Gamma_{m,n}}\\
            &- \inprod{\jump{\bus}}{\avg{\bn\cdot\nu\nabla\bu}}_{\Gamma_{m,n}}\\
            &+ \inprod{\gamma\Delta\bx^{-1}\jump{\bus}}{\jump{\bu}}_{\Gamma_{m,n}}
        \end{split}
    \end{equation}
    \begin{equation}
        \begin{split}
            \form{K_{m, di}}{\bus}{\bu} =&
            - \inprod{\bn\cdot\nu\nabla\bus}{\bu}_{\Gamma_{m,di}}\\
            &- \inprod{\bus}{\bn\cdot\nu\nabla\bu}_{\Gamma_{m,di}}\\
            &+ \inprod{\gamma\Delta\bx^{-1}\bus}{\bu}_{\Gamma_{m,di}},
        \end{split}
    \end{equation}
\end{subequations}
\todo{Shouldn't we use $\boldsymbol{g}_{di}$ here because Dirichlet boundary is involved?\\
KC: boundary condition is enforced weakly, so this term works as a sort of forcing term,
paired with a right-hand side with $\bg_{di}$.}
where the interface penalty strength $\gamma = \nu(s+1)^2$ is used for $\rU_{m,s}$
to ensure the stability of the solution~\cite{Toselli2002}, for both FOM and ROM.
The mixed bilinear form $B$ approximates the velocity divergence,
\begin{subequations}\label{eq:Bform}
    \begin{equation}
        \form{B_m}{\ps}{\bu} = -\inprod{\ps_m}{\nabla\cdot\bu_m}_{\Omega_m}
    \end{equation}
    \begin{equation}
        \form{B_{m,n}}{\ps}{\bu} = \inprod{\avg{\ps}}{\jump{\bn\cdot\bu}}_{\Gamma_{m,n}}
    \end{equation}
    \begin{equation}
        \form{B_{m, di}}{\ps}{\bu} = \inprod{\ps}{\bn\cdot\bu}_{\Gamma_{m,di}}.
    \end{equation}
\end{subequations}\todo{It is interesting that gradient of velocity is involved in the subdomain functionals, while gradient of pressure is involved in the interface functionals. Is some sort of integration by parts applied here?\\
KC: I'm not following this question so much, since there is no pressure gradient here. But interface functional does come from integration-by-part of pressure gradient and its symmetric interior penalty.}
\par
In this study, we chose the standard nonlinear form of the nonlinear advection~\cite{Elman2014},
\begin{equation}\label{eq:Cform}
    \form{C_m}{\bus}{\bu} = \inprod{\bus_m}{\bu_m\cdot\nabla\bu_m}_{\Omega_m},
\end{equation}
with $C_{m,n}=C_{m,di}=0$.
Note that the interface and boundary conditions are handled only with the linear operators in (\ref{eq:Kform}) and (\ref{eq:Bform}).
While this may cause numerical instability for high Reynolds numbers,
steady solutions for Navier-Stokes flow can \PR{physically} exist only up to modest Reynolds number~\cite{Elman2014}.
For flows past blunt bodies, steady flow stably exists only up to $\Re\lesssim40$ based on the object length scale~\cite{lienhard1966synopsis}.
In our study, the chosen formulation stably computes the solution up to $\mathrm{Re}\le50$ for our test cases.
As an alternative, a DG formulation with local Lax-Friedrichs flux can be used~\cite{fehn2018robust},
where $C_{m,n}$ and $C_{m,di}$ are non-zero.
\par
Just as for interface condition,
the boundary conditions are enforced weakly and
the right-hand side has the boundary linear forms,
\begin{subequations}
    \begin{equation}
        \form{L_u}{\bus}{\bg} = \inprod{\bus}{\bg_{ne}}_{\partial\Omega_{ne}} + \inprod{\bus}{\gamma\Delta\bx^{-1}\bg_{di}}_{\partial\Omega_{di}} + \inprod{\bn\cdot\nu\nabla\bus}{\bg_{di}}_{\partial\Omega_{di}}
    \end{equation}
    \begin{equation}
        \form{L_p}{\ps}{\bg} = \inprod{\ps}{\bn\cdot\bg_{di}}_{\partial\Omega_{di}}.
    \end{equation}
\end{subequations}
By having weakly-enforced boundary operators,
each boundary condition on each face of the domain can be selectively applied depending on problems,
without having to precompute all possible combinations.
In this study, no external forcing is considered; \PR{however, it can be readily extended to include non-zero forcing}.
\par
These functionals can also be written as matrix-vector inner products,
which explicitly show the operators that will be projected onto the reduced linear subspace in the subsequent section.
For the viscous flux operator,
\begin{subequations}\label{eq:vis-mat}
    \begin{equation}
        \form{K_m}{\bus}{\bu} = \bu_m^{\dagger\top}\bK_m\bu_m
    \end{equation}
    \begin{equation}
        \form{K_{m,n}}{\bus}{\bu} = 
        \begin{pmatrix}
            \bu_m^{\dagger\top} & \bu_n^{\dagger\top}
        \end{pmatrix}
        \begin{pmatrix}
            \bK_{mm} & \bK_{mn} \\
            \bK_{nm} & \bK_{nn} \\
        \end{pmatrix}
        \begin{pmatrix}
            \bu_m \\ \bu_n
        \end{pmatrix}
    \end{equation}
    \begin{equation}
        \form{K_{m,di}}{\bus}{\bu} = \bu_m^{\dagger\top}\bK_{m,di}\bu_m,
    \end{equation}
\end{subequations}
where $\bK_m, \bK_{m,di} \in \rR^{N_{u,m}\times N_{u,m}}$ with $N_{u,m}\equiv\dim(\bu_m)$ degrees of freedom for $\bu_m$, and $\bK_{ij} \in \rR^{N_{u,i}\times N_{u,j}}$.
For the velocity divergence operator,
\begin{subequations}\label{eq:div-mat}
    \begin{equation}\label{eq:div-comp}
        \form{B_m}{\ps}{\bu} = \bp_m^{\dagger\top}\bB_m\bu_m
    \end{equation}
    \begin{equation}
        \form{B_{m,n}}{\ps}{\bu} = 
        \begin{pmatrix}
            \bp_m^{\dagger\top} & \bp_n^{\dagger\top}
        \end{pmatrix}
        \begin{pmatrix}
            \bB_{mm} & \bB_{mn} \\
            \bB_{nm} & \bB_{nn} \\
        \end{pmatrix}
        \begin{pmatrix}
            \bu_m \\ \bu_n
        \end{pmatrix}
    \end{equation}
    \begin{equation}
        \form{B_{m,di}}{\ps}{\bu} = \bp_m^{\dagger\top}\bB_{m,di}\bu_m,
    \end{equation}
\end{subequations}
where $\bB_m, \bB_{m,di} \in \rR^{N_{p,m}\times N_{u,m}}$ with $N_{p,m}\equiv\dim(\bp_m)$ degrees of freedom for $\ps_m$, and $\bB_{ij} \in \rR^{N_{p,i}\times N_{u,j}}$.
The nonlinear advection operator can be similarly expressed as vector inner products,
\begin{equation}\label{eq:adv-mat}
    \form{C_m}{\bus}{\bu} = \bu_m^{\dagger\top}\bC[\bu_m]
\end{equation}
with $\bC_m: \rR^{N_{u,m}} \to \rR^{N_{u,m}}$.
For linear form right-hand side terms,
\begin{subequations}\label{eq:rhs-vec}
    \begin{equation}
        \form{L}{\bus}{\f} = \bu^{\dagger\top}\bL[\f]
    \end{equation}
    \begin{equation}
        \form{L_{u}}{\bus}{\f} = \bu^{\dagger\top}\bL_{u}[\bg]
    \end{equation}
    \begin{equation}
        \form{L_{p}}{\bqs}{\f} = \bp^{\dagger\top}\bL_{p}[\bg],
    \end{equation}
\end{subequations}
where $\bL, \bL_u: L_2(\Omega_m)^d\to\rR^{N_u}$ and $\bL_{p}: L_2(\Omega)^d\to\rR^{N_p}$
with $N_u\equiv\dim{\bu}$ and $N_p\equiv\dim{\bp}$ degrees of freedom for the global solution $\bu$ and $\bp$, respectively.
\par
The full-order model operators in (\ref{eq:vis-mat}–\ref{eq:rhs-vec}) establish the foundation
for creating reduced-order model components,
which can then be integrated into a global-scale reduced-order model.
The process for constructing and assembling the reduced-order model is detailed subsequently.

\subsection{Linear subspace approximation}
With all subdomains categorized with a few reference domains (\ref{eq:ref-domain}),
we seek to approximate the subdomain solutions $\bu_m$ and $\bp_m$
on a low-dimensional linear subspace for each reference domain $\bOmega_r$.
In other words, $\bu_m$ and $\bp_m$ are approximated as
linear combinations of given orthonormal basis vectors,
\begin{equation}\label{eq:linear-subspace}
    \bu_r \approx \bPhi_{u,r}\hbu_r
    \qquad\text{and}\qquad
    \bp_r \approx \bPhi_{p,r}\hbp_r,
\end{equation}
with $\bPhi_{u,r}\in\rR^{N_{u,r}\times R_{u,r}}$ and $\bPhi_{p,r}\in\rR^{N_{p,r}\times R_{p,r}}$
the basis matrix for velocity and pressure, respectively.
The coefficients for basis vectors $\hbu_r \in \rR^{R_{u,r}}$ and $\hbp_r \in \rR^{R_{p,r}}$
are then the ROM solution variable for the reference domain $r$.
$R_{u,r}$ and $R_{p,r}$ are the number of basis vectors on the reference domain $r$ for velocity and pressure, respectively.
For an effective reduced order model,
$\bu_r$ and $\bp_r$ must be well approximated with small $R_{u,r}\ll N_{u,r}$ and $R_{p,r} \ll N_{p,r}$,
without compromising the accuracy.
For parameterized problems, such bases can be well identified from snapshots in given reference domains
via proper orthogonal decomposition, which is described subsequently.
\par
The global-scale solution $\bu$ and $\bp$ are approximated in each subdomain,
\begin{subequations}\label{eq:global-rom-sol}
    \begin{equation}
        \bu \approx \bPhi_u{\hbu} \equiv
        \begin{pmatrix}
            \bPhi_{u,r(1)} & & & & \\
            & \cdots & & & \\
            & & \bPhi_{u,r(m)} & & \\
            & & & \cdots & \\
            & & & & \bPhi_{u,r(M)} \\
        \end{pmatrix}
        \left(\begin{array}{c}
            \hbu_1 \\
            \vdots \\
            \hbu_m \\
            \vdots \\
            \hbu_M \\
        \end{array}\right)
    \end{equation}    
    \begin{equation}
        \bp \approx \bPhi_p{\hbp} \equiv
        \begin{pmatrix}
            \bPhi_{p,r(1)} & & & & \\
            & \cdots & & & \\
            & & \bPhi_{p,r(m)} & & \\
            & & & \cdots & \\
            & & & & \bPhi_{p,r(M)} \\
        \end{pmatrix}
        \left(\begin{array}{c}
            \hbp_1 \\
            \vdots \\
            \hbp_m \\
            \vdots \\
            \hbp_M \\
        \end{array}\right),
    \end{equation}    
\end{subequations}
where $r(m)$ is the reference domain corresponding to the subdomain $m$.

\subsubsection{Proper orthogonal decomposition}


We consider $S_r$ sample velocity-pressure snapshots given for the reference domain $r$.
The snapshots can be obtained via solving (\ref{eq:dd-gov}) on the single reference domain $\bOmega_r$
or on a small-$M$ global domain that contains $\bOmega_r$.
Specific sampling procedures will be introduced
in Section~\ref{sec:results}.
All sample solutions for the decomposed system (\ref{eq:dd-gov}) are first
categorized by their reference domains, constituting the corresponding snapshot matrices
$\bU_r\in\rR^{N_{u,r}\times S_r}$ and $\bP_r\in\rR^{N_{p,r}\times S_r}$
\begin{subequations}\label{eq:snapshots}
    \begin{equation}
        \bU_r =
        \left(\begin{array}{c:c:c}
            & & \\
            {}_1\bu_r & {}_2\bu_r & \cdots\\
            & & \\
        \end{array}\right)
        \qquad
        \bP_r =
        \left(\begin{array}{c:c:c}
            & & \\
            {}_1\bp_r & {}_2\bp_r & \cdots\\
            & & \\
        \end{array}\right),
    \end{equation}
    where all sample solutions ${}_i\bu_r$ and ${}_i\bp_r$ lie on the reference domain $\bOmega_r$, i.e.,
    \begin{equation}
        \mathrm{dom}({}_i\bu_r) = \mathrm{dom}({}_i\bp_r) = \bOmega_r \qquad \forall i.
    \end{equation}
\end{subequations}
\par
The well-known proper orthogonal decomposition (POD) is used to decompose $\bU_r$ and $\bP_r$~\cite{berkooz1993proper,Chatterjee2000,Liang2002,Rowley2017},
\begin{equation}\label{eq:pod}
    \bU_r = \tilde{\bPhi}_{u,r}\tilde{\bS}_{u,r}\tilde{\bV}_{u,r}^{\top}
    \qquad
    \bP_r = \tilde{\bPhi}_{p,r}\tilde{\bS}_{p,r}\tilde{\bV}_{p,r}^{\top}.
\end{equation}
The orthonormal matrices $\tilde{\bPhi}_{u,r}\in\rR^{N_{u,r}\times S_r}$ and $\tilde{\bPhi}_{p,r}\in\rR^{N_{p,r}\times S_r}$
contain left singular vectors,
which can be used as basis vectors for the linear space of velocity and pressure, respectively.
The diagonal matrices $\tilde{\bS}_{u,r}, \tilde{\bS}_{p,r}\in\rR^{S_r\times S_r}$
contain singular values for corresponding singular vectors,
where each singular value indicates the total energy of the corresponding singular vector within the snapshot matrix.
In other words, the singular vectors with the larger singular values correspond to the most dominant physics modes representing the snapshots.
The matrices $\tilde{\bV}_{u,r}, \tilde{\bV}_{p,r}\in\rR^{S_{r}\times S_r}$ contain the right singular vectors,
\KC{which contain temporal behavior of the left singular vectors.
While these are often used in the context of space--time ROM formulation~\cite{choi2019space, choi2021space, kim2021efficient},
they do not play a major role here.}
\par
The singular vectors in $\tilde{\bPhi}_{u,r}$ and $\tilde{\bPhi}_{p,r}$ are conveniently ordered via POD in the descending order of their singular values,
thus according to their importance within the snapshot matrix.
This leads to the common choice of the ROM bases
being the first $R_{u,r}$ (or $R_{p,r}$) column vectors of $\tilde{\bPhi}_{u,r}$ (or $\tilde{\bPhi}_{p,r}$), i.e.
\begin{equation}\label{eq:rom-basis1}
    \bPhi_{u,r} = \trun_{R_{u,r}}(\tilde{\bPhi}_{u,r})
    \qquad
    \bPhi_{p,r} = \trun_{R_{p,r}}(\tilde{\bPhi}_{p,r}),
\end{equation}
where $\trun_R(\bPhi)$ is the truncated matrix with the first $R$ columns of $\bPhi$.
With the corresponding singular value matrix $\bS\in\rR^{S\times S}$, the expected accuracy with the first $R$ singular vectors
can be estimated with the missing energy ratio,
\begin{equation}\label{eq:pod-eps}
    \epsilon(\bS, R) = 1 - \frac{\sum_s^{R}{}_s\sigma}{\sum_s^{S}{}_s\sigma},
\end{equation}
where ${}_s\sigma$ is the $s$-th diagonal entry of $\bS$.
\par
However, for the steady incompressible Navier-Stokes equations, this simple choice can be problematic for a stable pressure prediction.
The indefinite nature of the system necessitates the solution space
to satisfy LBB condition~\cite{Babuvska1971,Brezzi1974,Ladyzhenskaya1963,Taylor1973}.
For the full-order model in Section~\ref{subsec:dg-dd}, this is satisfied via the choice of Taylor-Hood finite element space.
The same physics condition must be satisfied for the ROM basis space as well,
though the simple choice (\ref{eq:rom-basis1}) violates it:
the snapshots in $\bU_r$ are divergence-free as solutions of (\ref{eq:dd-gov-div}),
and the resulting bases $\bPhi_{u,r}$ from $\bU_r$ are likewise divergence-free.
In order to uniquely determine and stably compute the pressure,
the basis (\ref{eq:rom-basis1}) must be augmented with compressible components of the velocity field.
This augmentation procedure is described subsequently.


\subsubsection{Supremizer enrichment}\label{subsec:supreme}
We introduce the supremizer enrichment procedure in this section,
while referring readers to Ballarin \textit{et al.}~\cite{ballarin2015supremizer} for the theoretical background.
\par
We seek an augmented POD basis with an additional velocity basis $\bZ_r\in\rR^{N_{u,r}\times Z_r}$,
\begin{equation}\label{eq:vel-basis-sup}
    \bPhi_{u,r} =
    \begin{pmatrix}
        \trun_{R_{u,r}}(\tilde{\bPhi}_{u,r}) & \bZ_r
    \end{pmatrix}
    \in \rR^{N_{u,r}\times(R_{u,r} + Z_r)}.
\end{equation}
In order to obtain a compressible velocity component,
we compute the \textit{supremizer} of the pressure POD basis using the divergence operator (\ref{eq:div-comp}),
\begin{equation}\label{eq:supreme}
    \tilde{\bZ}_r = \bB_r^{\top}\trun_{Z_r}(\tilde{\bPhi}_{p,r}).
\end{equation}
Ballarin \textit{et al.}~\cite{ballarin2015supremizer} adopted a slightly different POD basis,
using generalized SVD with the mass matrix of velocity finite element space.
Their definition of supremizer is thus slightly different from (\ref{eq:supreme}),
where the mass matrix is multiplied on the left-hand side of (\ref{eq:supreme}).
In our study, this difference did not make a notable impact for a stable pressure prediction.
\par
In order to keep the orthonormality of the basis,
we perform the modified Gram-Schmidt procedure~\cite{cheney2009linear} with both POD basis and the supremizers,
\begin{equation}
    \begin{pmatrix}
        \trun_{R_{u,r}}(\tilde{\bPhi}_{u,r}) & \bZ_r
    \end{pmatrix}
    =
    \mathrm{GS}
    \begin{pmatrix}
        \trun_{R_{u,r}}(\tilde{\bPhi}_{u,r}) & \tilde{\bZ}_r
    \end{pmatrix}.
\end{equation}
Since $\tilde{\bPhi}_{u,r}$ is already orthonormal,
in practice only $\tilde{\bZ}_r$ is orthonormalized against $\tilde{\bPhi}_{u,r}$ and itself.
\todo{I am wondering if this Gram-Schmidt reveals that some of $\tilde{\boldsymbol{Z}}_r$ vectors are already in the subspace of $\boldsymbol{\Phi}_{u,r}$.\\
KC: Due to our choice of finite-element space and DG formulation, the solution is not exactly divergence-free, especially near boundaries.
So I think there are some compressible components in the POD basis already, but not sufficient to effectively stabilize our pressure basis.
}
\par
\KC{Ballarin \textit{et al.}~\cite{ballarin2015supremizer} proved and demonstrated
that each supremizer in (\ref{eq:supreme}) maximizes the stability factor of the corresponding pressure mode (thus supremizer).
This strongly suggests that at least the same number of supremizers should be used as pressure basis vectors, i.e. $Z_r \ge R_{p,r}$.}
While augmentation with $\bZ_r$ stabilizes the pressure,
the increase in basis size also increases the computational cost, reducing the effectiveness of ROM.
\KC{Therefore, following a rule of thumb suggested from numerical experiments by Ballarin \textit{et al.}~\cite{ballarin2015supremizer},
we adopt to use the minimum size of supremizers $Z_r = R_{p,r}$ throughout our study.} 

\subsection{POD-Galerkin ROM}

\subsubsection{Overview and linear operators}\label{subsubsec:galerkin-linear}
We project the decomposed system (\ref{eq:weak-fom-gov}) onto the linear subspaces of $\bPhi_u$ and $\bPhi_p$.
In particular, we consider Galerkin projection where the test function is also approximated with the same linear subspace as for (\ref{eq:global-rom-sol}),
\begin{equation}\label{eq:global-rom-test}
    \bus \approx \bPhi_u\hbus
    \qquad
    \bps \approx \bPhi_p\hbps.
\end{equation}
\par
The FOM operators (\ref{eq:vis-mat}-\ref{eq:rhs-vec}) are then reduced to
ROM operators for ROM variable $\hbu$ and $\hbp$.
For the viscous flux operator (\ref{eq:vis-mat}),
\begin{subequations}\label{eq:rom-vis-mat}
    \begin{equation}
        \begin{split}
            \form{K_m}{\bus}{\bu}
            &\approx \form{\hK_m}{\bus}{\bu}\\
            &\equiv \hbu_m^{\dagger\top}\hbK_m\hbu_m
            \equiv \hbu_m^{\dagger\top}(\bPhi_{u,r(m)}^{\top}\bK_m\bPhi_{u,r(m)})\hbu_m
        \end{split}
    \end{equation}
    \begin{equation}
    \begin{split}
        \form{K_{m,n}}{\bus}{\bu}
        &\approx \form{\hK_{m,n}}{\hbus}{\hbu}\\
        &\equiv 
        \begin{pmatrix}
            \hbu_m^{\dagger\top} & \hbu_n^{\dagger\top}
        \end{pmatrix}
        \begin{pmatrix}
            \hbK_{mm} & \hbK_{mn} \\
            \hbK_{nm} & \hbK_{nn} \\
        \end{pmatrix}
        \begin{pmatrix}
            \hbu_m \\ \hbu_n
        \end{pmatrix}\\
        &\equiv 
        \begin{pmatrix}
            \hbu_m^{\dagger\top} & \hbu_n^{\dagger\top}
        \end{pmatrix}
        \begin{pmatrix}
            \bPhi_{u,r(m)}^{\top}\bK_{mm}\bPhi_{u,r(m)} & \bPhi_{u,r(m)}^{\top}\bK_{mn}\bPhi_{r(n)} \\
            \bPhi_{u,r(n)}^{\top}\bK_{nm}\bPhi_{u,r(m)} & \bPhi_{u,r(n)}^{\top}\bK_{nn}\bPhi_{r(n)} \\
        \end{pmatrix}
        \begin{pmatrix}
            \hbu_m \\ \hbu_n
        \end{pmatrix}
    \end{split}
    \end{equation}
    \begin{equation}
        \begin{split}
            \form{K_{m,di}}{\bus}{\bu}
            &\approx \form{\hK_{m,di}}{\hbus}{\hbu}\\
            &\equiv \hbu_m^{\dagger\top}\hbK_{m,di}\hbu_m
            \equiv \hbu_m^{\dagger\top}\left(\bPhi_{u,r(m)}^{\top}\bK_{m,di}\bPhi_{u,r(m)}\right)\hbu_m.
        \end{split}
    \end{equation}
\end{subequations}
For the velocity divergence operator (\ref{eq:div-mat}),
\begin{subequations}\label{eq:rom-div-mat}
    \begin{equation}\label{eq:comp-bilinear-mat-reduced:Bm}
    \begin{split}
        \form{B_m}{\bps}{\bu}
        &\approx \form{\hB_m}{\hbps}{\hbu}\\
        &\equiv \hbp_m^{\dagger\top}\hbB_m\hbu_m
        \equiv \hbp_m^{\dagger\top}\left(\bPhi_{p,r(m)}^{\top}\bB_m\bPhi_{u,r(m)}\right)\hbu_m
    \end{split}
    \end{equation}
    \begin{equation}\label{eq:comp-bilinear-mat-reduced:Bmn}
    \begin{split}
        \form{B_{m,n}}{\bps}{\bu}
        &\approx \form{\hB_{m,n}}{\hbps}{\hbu}\\
        &\equiv 
        \begin{pmatrix}
            \hbp_m^{\dagger\top} & \hbp_n^{\dagger\top}
        \end{pmatrix}
        \begin{pmatrix}
            \hbB_{mm} & \hbB_{mn} \\
            \hbB_{nm} & \hbB_{nn} \\
        \end{pmatrix}
        \begin{pmatrix}
            \hbu_m \\ \hbu_n
        \end{pmatrix}\\
        &\equiv 
        \begin{pmatrix}
            \hbp_m^{\dagger\top} & \hbp_n^{\dagger\top}
        \end{pmatrix}
        \begin{pmatrix}
            \bPhi_{p,r(m)}^{\top}\bB_{mm}\bPhi_{u,r(m)} & \bPhi_{p,r(m)}^{\top}\bB_{mn}\bPhi_{u,r(n)} \\
            \bPhi_{p,r(n)}^{\top}\bB_{nm}\bPhi_{u,r(m)} & \bPhi_{p,r(n)}^{\top}\bB_{nn}\bPhi_{u,r(n)} \\
        \end{pmatrix}
        \begin{pmatrix}
            \hbu_m \\ \hbu_n
        \end{pmatrix}
    \end{split}
    \end{equation}
    \begin{equation}
        \begin{split}
            \form{B_{m,di}}{\bps}{\bu}
            &\approx \form{\hB_{m,di}}{\hbps}{\hbu}\\
            &\equiv \hbp_m^{\dagger\top}\hbB_{m,di}\hbu_m
            \equiv \hbp_m^{\dagger\top}\left(\bPhi_{p,r(m)}^{\top}\bB_{m,di}\bPhi_{u,r(m)}\right)\hbu_m.
        \end{split}
        \end{equation}
\end{subequations}
Likewise, the vectors of the linear form in (\ref{eq:rhs-vec}) are also reduced to
\begin{subequations}\label{eq:rom-rhs-vec}
    \begin{equation}
    \begin{split}
        \form{L}{\bus}{\f}
        &\approx \form{\hL}{\hbus}{\f}\\
        &\equiv \hbu^{\dagger\top}\hbL[\f]
        \equiv \hbu^{\dagger\top}\left(\bPhi_{u}^{\top}\bL[\f]\right)
    \end{split}
    \end{equation}
    \begin{equation}
    \begin{split}
        \form{L_u}{\bus}{\f}
        &\approx \form{\hL_u}{\hbus}{\bg}\\
        &\equiv \hbu^{\dagger\top}\hbL_u[\bg]
        \equiv \hbu^{\dagger\top}\left(\bPhi_u^{\top}\bL_u[\bg]\right)
    \end{split}
    \end{equation}
    \begin{equation}
    \begin{split}
        \form{L_p}{\bps}{\f}
        &\approx \form{\hL_p}{\hbps}{\f}\\
        &\equiv \hbp^{\dagger\top}\hbL_p[\f]
        \equiv \hbp^{\dagger\top}\left(\bPhi_p^{\top}\bL_p[\bg]\right).
    \end{split}
    \end{equation}
\end{subequations}

\subsubsection{Nonlinear advection operator: (1) tensorial approach}
With (\ref{eq:global-rom-sol}),
the nonlinear advection operator (\ref{eq:adv-mat}) is approximated as
\begin{equation}\label{eq:nladvec-rom1}
    \form{C_m}{\bus}{\bu}
    \approx
    \hbu_m^{\dagger\top}\bPhi_{u,r(m)}^{\top}\bC[\bPhi_{u,r(m)}\hbu_m].
\end{equation}
Unlike linear operators in Section~\ref{subsubsec:galerkin-linear},
this advection ROM operator cannot be readily assembled before prediction due to its nonlinearity.
Naive computation of (\ref{eq:nladvec-rom1}) would involve
evaluation of the operator $\bC[\cdot]$ at every element and quadrature point,
which does not reduce the computational cost at all.
\par
An alternative to achieve speed-up similar to the linear operators is the tensorial approach.
We expand $\bu$ and $\bus$ in the advection operator (\ref{eq:Cform}) with (\ref{eq:global-rom-sol}) and (\ref{eq:global-rom-test}),
\begin{equation}\label{eq:rom-adv-tensor}
\begin{split}
    \form{C_m}{\bus}{\bu} &\approx
    \sum_{i,j,k}^{R_{u,r(m)}}
    \hus_i\hu_j\hu_k
    \inprod{{}_i\bphi_{u,r(m)}}{{}_j\bphi_{u,r(m)}\cdot\nabla {}_k\bphi_{u,r(m)}}_{\Omega_m}\\
    &\equiv
    \hbu^{\dagger\top}\uhbC_{r(m)}:\hbu\hbu
    \equiv \form{\hC_{m}}{\hbus}{\hbu},
\end{split}
\end{equation}
where the entries of the third-order tensor $\uhbC_r\in\rR^{R_{u,r}\times R_{u,r}\times R_{u,r}}$ are combinations of three velocity basis vectors,
\begin{equation}
    \left(\uhbC_{r}\right)_{ijk} = \inprod{{}_i\bphi_{u,r(m)}}{{}_j\bphi_{u,r(m)}\cdot\nabla {}_k\bphi_{u,r(m)}}_{\Omega_m}.
\end{equation}
This way, the ROM operator for the nonlinear advection can be precomputed as for its linear counterparts.
In online prediction, the complexity of (\ref{eq:rom-adv-tensor}) scales as $\cO(R_{u,r}^3)$,
much faster compared to $\cO(R_{u,r}^2)$ for (\ref{eq:rom-vis-mat}) or (\ref{eq:rom-div-mat}).
However, we can still expect a significant speed-up as long as a moderate basis size $R_{u,r}$ is used.
\par
The tensorial approach exploits the fact that the advection operator (\ref{eq:Cform}) is quadratic with respect to $\bu$.
As a result, it is not extensible to general nonlinear terms.
In such cases, the empirical quadrature procedure is another alternative to reduce the computation cost of the nonlinear operator.
This is introduced subsequently.

\subsubsection{Nonlinear advection operator: (2) empirical quadrature procedure}

Here we employ another alternative to reduce the computational cost of the nonlinear operator,
following the work of Chapman \textit{et al.}~\cite{chapman2017accelerated}.
\par
In essence, EQP attempts to approximate integral quantities---particularly functionals such as (\ref{eq:Cform})---with a reduced numerical quadrature,
\begin{equation}
    \begin{split}
        \inprod{\bus}{\bu\cdot\nabla\bu}_{\Omega_m}
        &\approx
        \sum_{q=1}^{N_q}w_q\bus(\bx_q)^{\top}\left(\bu(\bx_q)\cdot\nabla\bu(\bx_q)\right)\\
        &\equiv
        \inprod{\bus}{\bu\cdot\nabla\bu}_{\Omega_m,EQP}.
    \end{split}
\end{equation}
where the empirical quadrature points $\{\bx_q\}_{q=1}^{N_q}$ are sub-sampled
from the quadrature points for FOM finite-element spaces $\rU_{s+1}$ and $\rP_{s}$,
with their weights $\{w_q\}_{q=1}^{N_q}$ calibrated accordingly.
For the modest accuracy expected from the reduced subspace approximation (\ref{eq:linear-subspace}),
much fewer quadrature points can be sufficient.
\par
We seek to find such empirical quadrature points and weights for the functional (\ref{eq:Cform})
at the reference domain $r$ with the given snapshots $\bU_r$ (\ref{eq:snapshots}) and the basis $\bPhi_{u,r}$ (\ref{eq:vel-basis-sup}).
Specifically, for $\cQ_r$ the set of all quadrature points in the reference domain $\bOmega_r$,
we seek a weight vector $\bw_r\in\rR^{n(\cQ_r)}$
whose element $w_{r,q}$ corresponds to the $q$-th quadrature point $\bx_q\in\cQ$, such that
\begin{subequations}\label{eq:eqp}
    \begin{equation}
        \bw_r = \arg\min\Vert\bw'_r\Vert_0,
    \end{equation}
    where the weight vector $\bw'_r$ satisfies
    \begin{equation}
        w'_{r,q} \ge 0 \qquad \forall q\in[1, n(\cQ_r)]
    \end{equation}
    \begin{equation}\label{eq:eqp-error}
        \begin{split}
            \sum_{b=1}^{R_{u,r}}\sum_{s=1}^{S_r} \bigg[\inprod{{}_b\bphi_{u,r}}{{}_s\bu_r\cdot\nabla{}_s\bu_r}_{\Omega_r}
                &- \inprod{{}_b\bphi_{u,r}}{{}_s\bu_r\cdot\nabla{}_s\bu_r}_{\Omega_r,EQP}\bigg]^2\\
            &\le \epsilon_{EQP}^2 \sum_{b=1}^{R_{u,r}}\sum_{s=1}^{S_r} \left[\inprod{{}_b\bphi_{u,r}}{{}_s\bu_r\cdot\nabla{}_s\bu_r}_{\Omega_r}\right]^2,
        \end{split}
    \end{equation}
    with ${}_b\bphi_{u,r}$ the $b$-th column vector of $\bPhi_{u,r}$ and ${}_s\bu_r$ the $s$-th column vector of $\bU_r$.
\end{subequations}
The minimization problem (\ref{eq:eqp}) is solved via the non-negative least-squares active-set algorithm~\cite{chapman2017accelerated, lawson1995solving}.
\par
By employing the reduced numerical quadrature,
EQP applies an additional layer of approximation compared to the tensorial approach.
Its accuracy is controlled by the condition (\ref{eq:eqp-error}) with a specified relative error threshold $\epsilon_{EQP}$.
While a proper value for $\epsilon_{EQP}$ is problem-dependent,
a rule of thumb is to match with the missing energy ratio (\ref{eq:pod-eps}) of the ROM basis,
\begin{equation}\label{eq:eqp-eps}
    \epsilon_{EQP} = \epsilon(\tilde{\bS}_{u,r}, R_{u,r}).
\end{equation}
In our demonstration, with this error criterion,
EQP does not compromise the accuracy compared to the tensorial approach,
while achieving a slightly cheaper computational cost.
\par
With the calibrated EQP points and weights,
the ROM nonlinear advection operator (\ref{eq:nladvec-rom1}) is approximated as
\begin{equation}\label{eq:rom-adv-eqp}
    \begin{split}
        \form{C_m}{\bus}{\bu}
        &\approx
        \form{\hC_m}{\hbus}{\hbu}\\
        &\equiv
        \sum_{i}^{R_{u,r(m)}}\hus_i\inprod{{}_i\bphi_{u,r(m)}}{\bu_{EQP,r}\cdot\nabla \bu_{EQP,r}}_{\Omega_m,EQP},
    \end{split}
\end{equation}
where $\bu_{EQP,r}\in\rR^{N_{r,q}}$ is the sampled FOM velocity vector on EQP points,
\begin{equation}
    \left(\bu_{EQP,r}\right)_q = \sum_i^{R_{u,r}}{}_i\bphi_{u,r}(\bx_q)\hu_i.
\end{equation}

\subsubsection{Global ROM system}
With (\ref{eq:rom-vis-mat}-\ref{eq:rom-rhs-vec}), (\ref{eq:rom-adv-tensor}) and (\ref{eq:rom-adv-eqp}) as ROM components,
the global-scale ROM is assembled analogous to (\ref{eq:weak-fom-gov}),
\begin{subequations}\label{eq:rom-gov}
    \begin{equation}
        \form{\hK}{\hbus}{\hbu} + \form{\hB}{\hbp}{\hbus} + \form{\hC}{\hbus}{\hbu}
        =
        \form{\hL}{\hbus}{\f} + \form{\hL_u}{\hbus}{\bg}
    \end{equation}
    \begin{equation}
        \form{\hB}{\hbps}{\hbu}
        =
        \form{\hL_p}{\hbps}{\bg},
    \end{equation}
\end{subequations}
where the assembled global scale ROM operators analogous to (\ref{eq:fom-global-op}) are
\begin{subequations}\label{eq:rom-global-op}
    \begin{equation}
        \form{\hK}{\hbus}{\hbu} =
        \sum\limits_{m=1}^{M}\form{\hK_m}{\hbus}{\hbu} + 
        \sum\limits_{(m,n)\in \rI}\form{\hK_{m,n}}{\hbus}{\hbu} +
        \sum\limits_{m=1}^{M}\form{\hK_{m,di}}{\hbus}{\hbu}
    \end{equation}
    \begin{equation}
        \form{\hB}{\hbps}{\hbu} =
        \sum\limits_{m=1}^{M}\form{\hB_m}{\hbps}{\hbu} + 
        \sum\limits_{(m,n)\in \rI}\form{\hB_{m,n}}{\hbps}{\hbu} +
        \sum\limits_{m=1}^{M}\form{\hB_{m,di}}{\hbps}{\hbu}
    \end{equation}
    \begin{equation}\label{eq:rom-global-adv}
        \form{\hC}{\hbus}{\hbu} = 
        \sum\limits_{m=1}^{M} \form{\hC_m}{\hbus}{\hbu}.
    \end{equation}
\end{subequations}
For (\ref{eq:rom-global-adv}),
either tensorial approach (\ref{eq:rom-adv-tensor}) or EQP (\ref{eq:rom-adv-eqp}) can be adopted.
The performance of the two approaches will be compared in Section~\ref{sec:results}.
The global ROM equation (\ref{eq:rom-gov}) is solved for $\hbu$ and $\hbp$ with respect to all $\hbus$ and $\hbp$.

\section{Results}\label{sec:results}

\subsection{Sample generation and POD bases}\label{subsec:bases}
The proposed method is demonstrated on flow past objects devised by Chung \textit{et al.}~\cite{chung2024train}.
The domain $\Omega=[0, L]^2$ is composed of 5 different unit-square reference domains $\rC=\{\bOmega_r\}_{r=1}^5$,
where each reference domain contains a distinct object in it: empty, circle, triangle, square and star.
The reference domain meshes from Chung \textit{et al.}~\cite{chung2024train} are used with a uniform refinement,
in order to resolve flow solutions at $\Re=25$.
The inflow velocity is set to be a constant velocity added with a sinusoidal perturbation,
\KC{
\begin{equation}\label{eq:inflow-bc}
    \bg_{di} = (g_1 + \Delta g_1\sin 2\pi(\mathbf{k}_1\cdot\bx + \theta_1), g_2 + \Delta g_2\sin 2\pi(\mathbf{k}_2\cdot\bx + \theta_2)).
\end{equation}
}
Based on signs of the inflow, Dirichlet boundary condition is set on upwind sides of the domain.
The opposite sides are set to be homogeneous Neumann condition,
\KC{as the downstream flow may not be well characterized with an analytic function.
At the surface of the object inside the domain, no-slip wall boundary $\bg_{di} = \mathbf{0}$ is set.
Sample snapshots are obtained from 2000 random reference domains of $2\times2$ arrays with each array having randomly chosen inflow velocity from a uniform distribution,
\begin{subequations}\label{eq:train-sample}
    \begin{equation}
        g_1, g_2 \in U[-1, 1]
    \end{equation}
    \begin{equation}
        \Delta g_1, \Delta g_2 \in U[-0.1, 0.1]
    \end{equation}
    \begin{equation}
        \mathbf{k}_1, \mathbf{k}_2 \in U[-0.5, 0.5]^2
    \end{equation}
    \begin{equation}
        \theta_1, \theta_2 \in U[0, 1].
    \end{equation}
\end{subequations}
}
The numerical experiments in this study are performed on an Intel Sapphire Rapids CPU, running on a single thread with 2GB memory and 2GHz clock speed~\cite{dane}.
\par
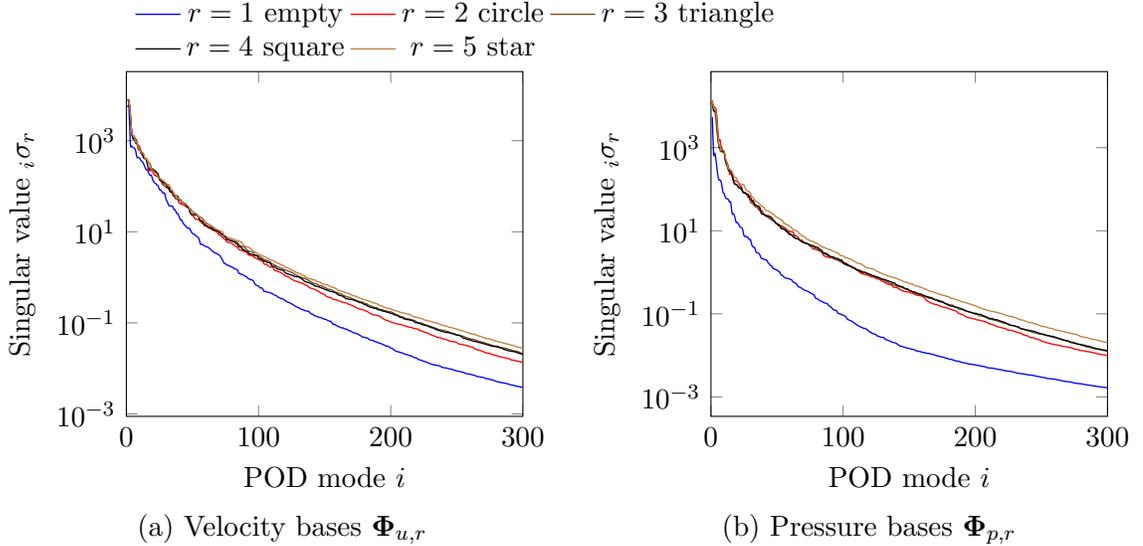
\begin{figure}[tbph]
    \begin{tikzpicture}[font=\small,]
    \begin{groupplot}[
        group style={
            group name = my plots,
            group size= 2 by 2,
            xlabels at =edge bottom,
            horizontal sep=2.5cm,
            vertical sep=2.2cm,
        },
        name=chung,
    ]    

        \nextgroupplot[
            height = 0.45\textwidth,
            width = 0.5\textwidth,
            ylabel={Singular value ${}_i\sigma_r$},
            xlabel={POD mode $i$},
            tick scale binop ={\times},
            ymode=log,
            xmin=0, xmax=300,
            legend style={
                draw=none, fill=none,
                at={(rel axis cs: 0., 1.0)},
                anchor=south west,
                legend columns=3,
            },
        ]
        
            \addplot+ [
                line width=0.5,
                smooth, solid,
                mark=none,
                select coords between index={0}{300},
            ]
            table [x expr=\coordindex+1, y index=0]{./data/ns_basis_empty_vel_sv.txt};
            \addplot+ [
                line width=0.5,
                smooth, solid,
                mark=none,
                select coords between index={0}{300},
            ]
            table [x expr=\coordindex+1, y index=0]{./data/ns_basis_square-circle_vel_sv.txt};
            \addplot+ [
                line width=0.5,
                smooth, solid,
                mark=none,
                select coords between index={0}{300},
            ]
            table [x expr=\coordindex+1, y index=0]{./data/ns_basis_square-triangle_vel_sv.txt};
            \addplot+ [
                line width=0.5,
                smooth, solid,
                mark=none,
                select coords between index={0}{300},
            ]
            table [x expr=\coordindex+1, y index=0]{./data/ns_basis_square-square_vel_sv.txt};
            \addplot+ [
                line width=0.5,
                brown,
                smooth, solid,
                mark=none,
                select coords between index={0}{300},
            ]
            table [x expr=\coordindex+1, y index=0]{./data/ns_basis_square-star_vel_sv.txt};

            \legend{$r=1$ empty, $r=2$ circle, $r=3$ triangle, $r=4$ square, $r=5$ star}

        \nextgroupplot[
            height = 0.45\textwidth,
            width = 0.5\textwidth,
            ylabel={Singular value ${}_i\sigma_r$},
            xlabel={POD mode $i$},
            tick scale binop ={\times},
            ymode=log,
            xmin=0, xmax=300,
        ]
        
            \addplot+ [
                line width=0.5,
                smooth, solid,
                mark=none,
                select coords between index={0}{300},
            ]
            table [x expr=\coordindex+1, y index=0]{./data/ns_basis_empty_pres_sv.txt};
            \addplot+ [
                line width=0.5,
                smooth, solid,
                mark=none,
                select coords between index={0}{300},
            ]
            table [x expr=\coordindex+1, y index=0]{./data/ns_basis_square-circle_pres_sv.txt};
            \addplot+ [
                line width=0.5,
                smooth, solid,
                mark=none,
                select coords between index={0}{300},
            ]
            table [x expr=\coordindex+1, y index=0]{./data/ns_basis_square-triangle_pres_sv.txt};
            \addplot+ [
                line width=0.5,
                smooth, solid,
                mark=none,
            ]
            table [x expr=\coordindex+1, y index=0]{./data/ns_basis_square-square_pres_sv.txt};
            \addplot+ [
                line width=0.5,
                brown,
                smooth, solid,
                mark=none,
                select coords between index={0}{300},
            ]
            table [x expr=\coordindex+1, y index=0]{./data/ns_basis_square-star_pres_sv.txt};
 
  \end{groupplot}
\node[below = 1.5cm of my plots c1r1.south west,
    anchor=west,
] {(a) Velocity bases $\bPhi_{u,r}$};
\node[below = 1.5cm of my plots c2r1.south west,
    anchor=west,
] {(b) Pressure bases $\bPhi_{p,r}$};
\end{tikzpicture}
%
    \caption{Singular value spectrum of the component POD bases: (a) velocity; and (b) pressure.}
    \label{fig:sv}
\end{figure}
Figure~\ref{fig:sv} shows the singular value spectra of the component POD bases.
For both velocity and pressure,
the spectra decays quickly with the POD mode,
indicating that the solutions can be effectively well represented on a low-dimensional subspace.
\begin{figure}[tbph]
    \begin{tikzpicture}[font=\small,]
    \begin{groupplot}[
        group style={
            group name = my plots,
            group size= 2 by 2,
            xlabels at =edge bottom,
            horizontal sep=2.5cm,
            vertical sep=2.2cm,
        },
        name=chung,
    ]    

        \nextgroupplot[
            height = 0.45\textwidth,
            width = 0.5\textwidth,
            ylabel={$\epsilon(\bS_{u,r}, R_{u,r})$},
            xlabel={$R_{u,r}$},
            tick scale binop ={\times},
            ymode=log,
            ytick={3e-3, 5e-3, 7e-3, 1e-2, 3e-2, 5e-2},
            yticklabels={$3\times10^{-3}$, $5\times10^{-3}$, $7\times10^{-3}$, $10^{-2}$, $3\times10^{-2}$, $5\times10^{-2}$},
            legend style={
                draw=none, fill=none,
                at={(rel axis cs: 0., 1.0)},
                anchor=south west,
                legend columns=3,
            },
        ]
        
            \addplot+ [
                line width=1.0,
                smooth, solid,
                mark=*,
                mark options={fill=white,solid},
            ]
            table [x index=0, y index=1]{./data/energy_ratio.txt};
            \addplot+ [
                line width=1.0,
                smooth, solid,
                mark=*,
                mark options={fill=white,solid},
            ]
            table [x index=0, y index=2]{./data/energy_ratio.txt};
            \addplot+ [
                line width=1.0,
                smooth, solid,
                mark=*,
                mark options={fill=white,solid},
            ]
            table [x index=0, y index=5]{./data/energy_ratio.txt};
            \addplot+ [
                line width=1.0,
                smooth, solid,
                mark=*,
                mark options={fill=white,solid},
            ]
            table [x index=0, y index=3]{./data/energy_ratio.txt};
            \addplot+ [
                line width=1.0,
                smooth, solid,
                mark=*,
                mark options={fill=white,solid},
                brown,
            ]
            table [x index=0, y index=4]{./data/energy_ratio.txt};

            \legend{$r=1$ empty, $r=2$ circle, $r=3$ triangle, $r=4$ square, $r=5$ star}

        \nextgroupplot[
            height = 0.45\textwidth,
            width = 0.5\textwidth,
            ylabel={$\epsilon(\bS_{p,r}, R_{p,r})$},
            xlabel={$R_{p,r}$},
            tick scale binop ={\times},
            ymode=log,
        ]
        
            \addplot+ [
                line width=1.0,
                smooth, solid,
                mark=*,
                mark options={fill=white,solid},
            ]
            table [x index=0, y index=6]{./data/energy_ratio.txt};
            \addplot+ [
                line width=1.0,
                smooth, solid,
                mark=*,
                mark options={fill=white,solid},
            ]
            table [x index=0, y index=7]{./data/energy_ratio.txt};
            \addplot+ [
                line width=1.0,
                smooth, solid,
                mark=*,
                mark options={fill=white,solid},
            ]
            table [x index=0, y index=10]{./data/energy_ratio.txt};
            \addplot+ [
                line width=1.0,
                smooth, solid,
                mark=*,
                mark options={fill=white,solid},
            ]
            table [x index=0, y index=8]{./data/energy_ratio.txt};
            \addplot+ [
                line width=1.0,
                smooth, solid,
                mark=*,
                mark options={fill=white,solid},
                brown,
            ]
            table [x index=0, y index=9]{./data/energy_ratio.txt};
 
  \end{groupplot}
\node[below = 1.5cm of my plots c1r1.south west,
    anchor=west,
] {(a) Velocity bases $\bPhi_{u,r}$};
\node[below = 1.5cm of my plots c2r1.south west,
    anchor=west,
] {(b) Pressure bases $\bPhi_{p,r}$};
\end{tikzpicture}
%
    \caption{Energy missing ratio (\ref{eq:pod-eps}) of the component POD bases: (a) velocity; and (b) pressure.}
    \label{fig:energy-ratio}
\end{figure}
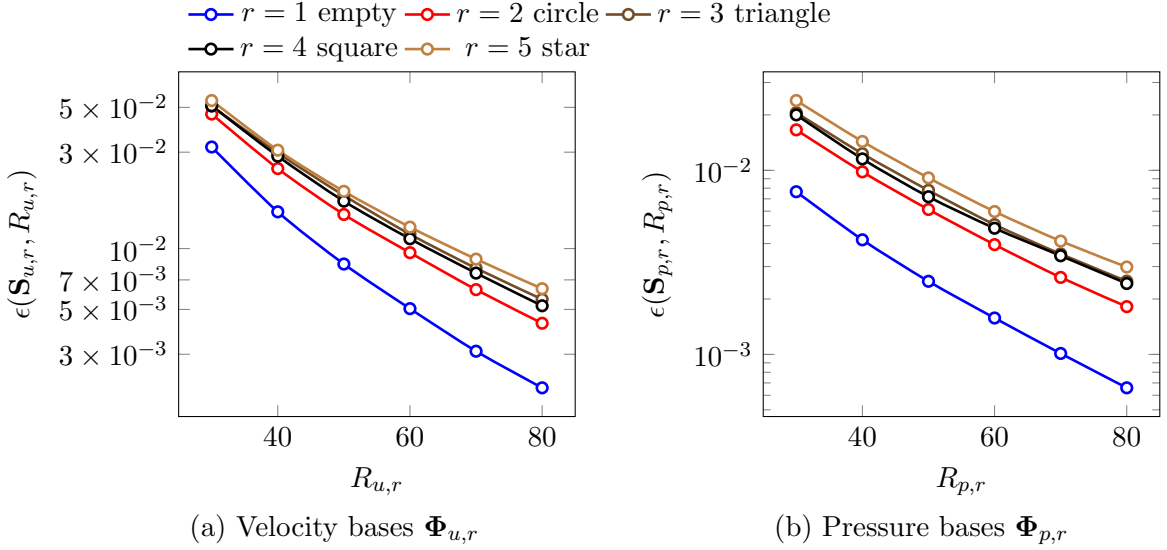
Figure~\ref{fig:energy-ratio}
shows the energy missing ratio (\ref{eq:pod-eps}) of the component POD bases.
Overall, pressure bases require fewer basis vectors than velocity to achieve the same energy missing ratio.
With $R_{u,r}=40$, velocity bases only miss less than $3\%$ for all reference components.

\subsection{Prediction performance over scale}
\begin{figure}[tbph]
    \begin{tikzpicture}[font=\small, spy using outlines={circle,black,magnification=6,size=1.5cm, connect spies}]
\pgfplotsset{set layers=standard}
    \begin{groupplot}[
        group style={
            group name = my plots,
            group size= 2 by 2,
            xlabels at =edge bottom,
            horizontal sep=3.5cm,
            vertical sep=2.5cm,
        },
        name=chung,
    ]    

        \nextgroupplot[
            height = 0.45\textwidth,
            width = 0.45\textwidth,
            ylabel={$x_2$},
            xlabel={$x_1$},
            tick scale binop ={\times},
            xmin = 0, xmax = 16,
            ymin = 0, ymax = 16,
            point meta min=0.0, point meta max=5.4,
            colorbar, colormap/jet,
            colorbar style={
                font=\scriptsize,
                xticklabel pos=upper,
                scaled y ticks=false,
                /pgf/number format/precision=4,
                at={(rel axis cs: 1.01, 0.)}, anchor=south west,
                xlabel=$\vert \bu\vert(\bx)$,
            }
        ]
        
            \edef\imagepath{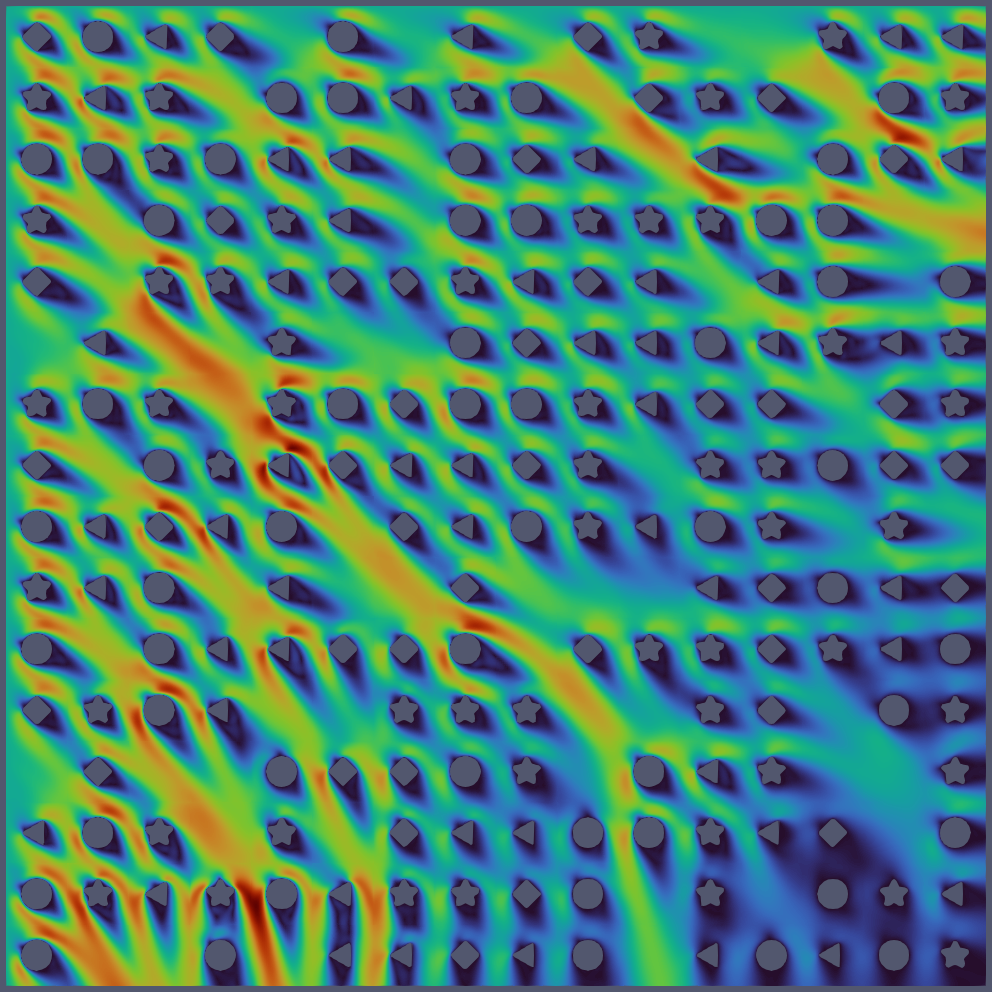}
            \addplot graphics[xmin=-0.1,xmax=16.1,ymin=-0.1,ymax=16.1]{\imagepath};


        \nextgroupplot[
            height = 0.45\textwidth,
            width = 0.45\textwidth,
            ylabel={$x_2$},
            xlabel={$x_1$},
            tick scale binop ={\times},
            xmin = 0, xmax = 16,
            ymin = 0, ymax = 16,
            point meta min=0.0, point meta max=0.43,
            colorbar, colormap/jet,
            colorbar style={
                font=\scriptsize,
                xticklabel pos=upper,
                scaled y ticks=false,
                /pgf/number format/precision=4,
                at={(rel axis cs: 2.8, 0.)}, anchor=south west,
                xlabel=$\epsilon_{\bu}(\bx)$,
            }
        ]
        
            \edef\imagepath{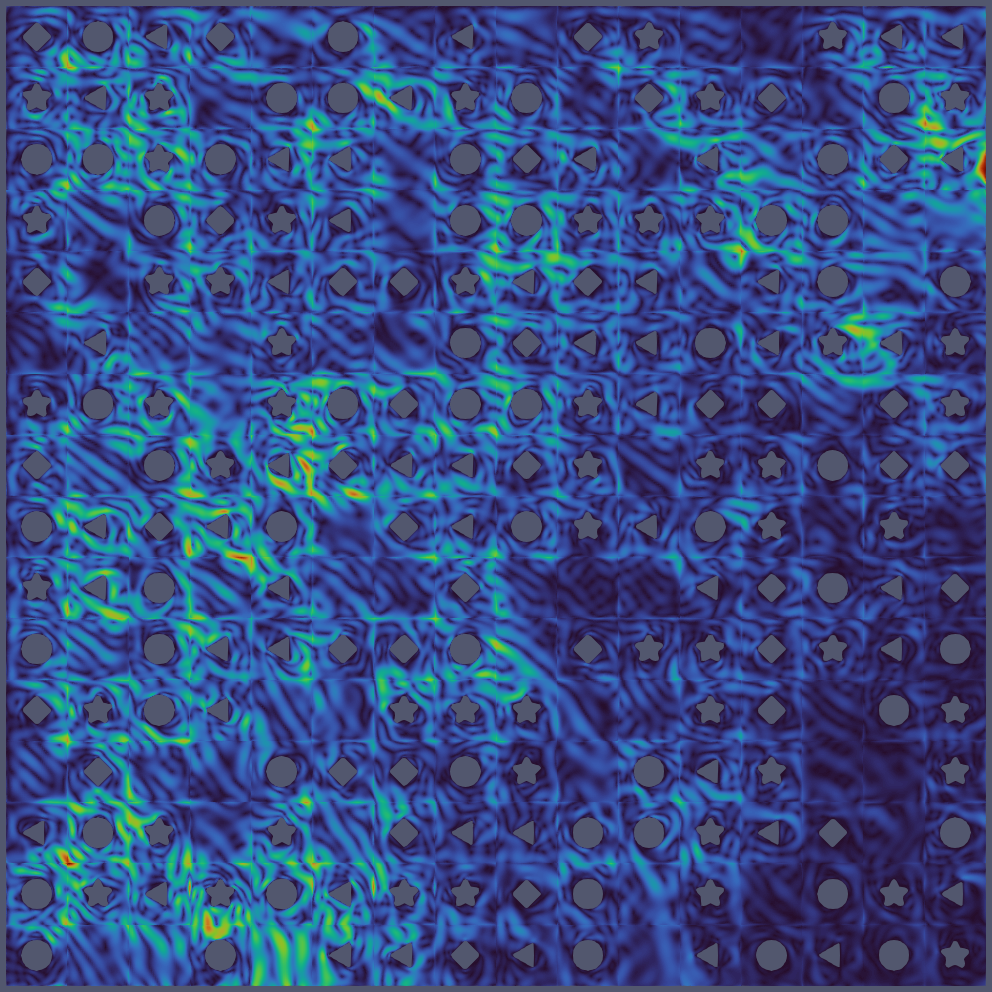}
            \addplot graphics[xmin=-0.1,xmax=16.1,ymin=-0.1,ymax=16.1]{\imagepath};

        \nextgroupplot[
            height = 0.45\textwidth,
            width = 0.45\textwidth,
            ylabel={$x_2$},
            xlabel={$x_1$},
            tick scale binop ={\times},
            xmin = 0, xmax = 16,
            ymin = 0, ymax = 16,
            point meta min=0.0, point meta max=52,
            colorbar, colormap/jet,
            colorbar style={
                font=\scriptsize,
                xticklabel pos=upper,
                scaled y ticks=false,
                /pgf/number format/precision=4,
                at={(rel axis cs: 1.01, -1.55)}, anchor=south west,
                xlabel=$p(\bx)$,
            }
        ]
        
            \edef\imagepath{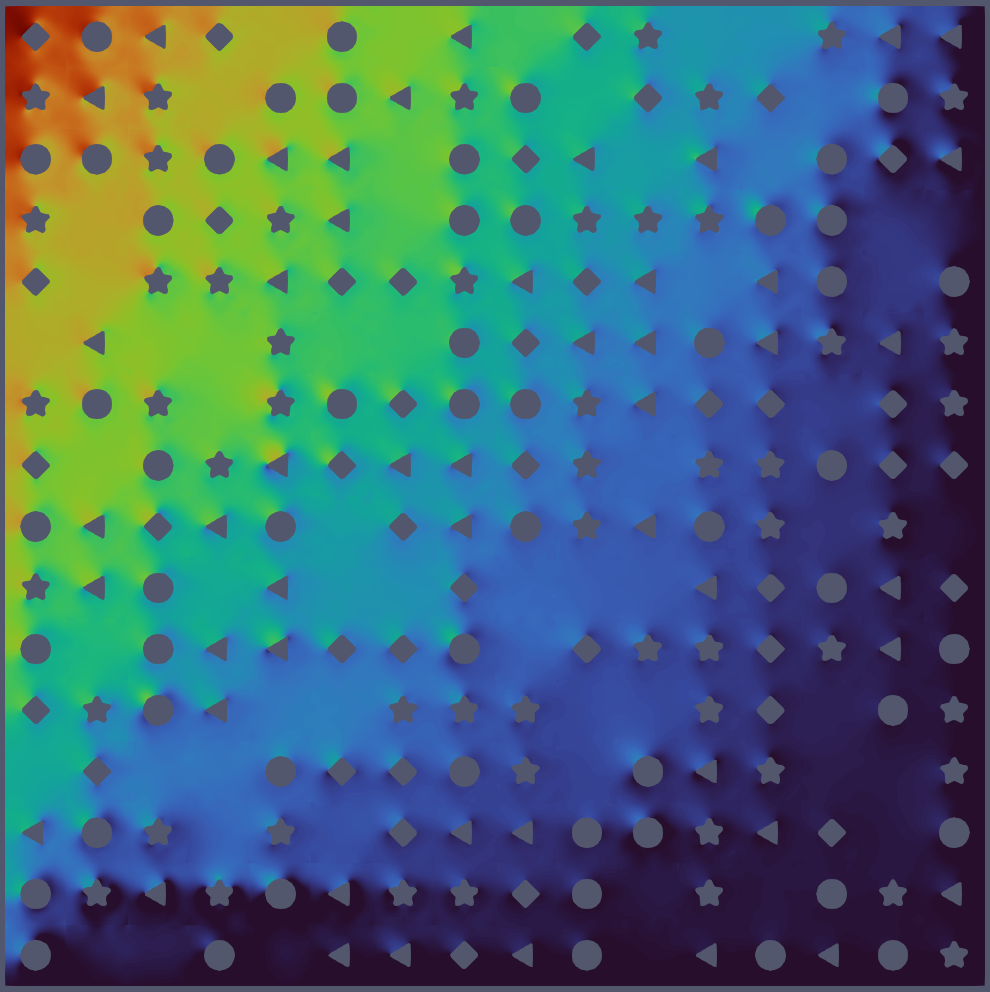}
            \addplot graphics[xmin=-0.1,xmax=16.1,ymin=-0.1,ymax=16.1]{\imagepath};

        \nextgroupplot[
            height = 0.45\textwidth,
            width = 0.45\textwidth,
            ylabel={$x_2$},
            xlabel={$x_1$},
            tick scale binop ={\times},
            xmin = 0, xmax = 16,
            ymin = 0, ymax = 16,
            point meta min=0.0, point meta max=2.5,
            colorbar, colormap/jet,
            colorbar style={
                font=\scriptsize,
                xticklabel pos=upper,
                scaled y ticks=false,
                /pgf/number format/precision=4,
                at={(rel axis cs: 2.8, -1.55)}, anchor=south west,
                xlabel=$\epsilon_p(\bx)$,
            }
        ]
        
            \edef\imagepath{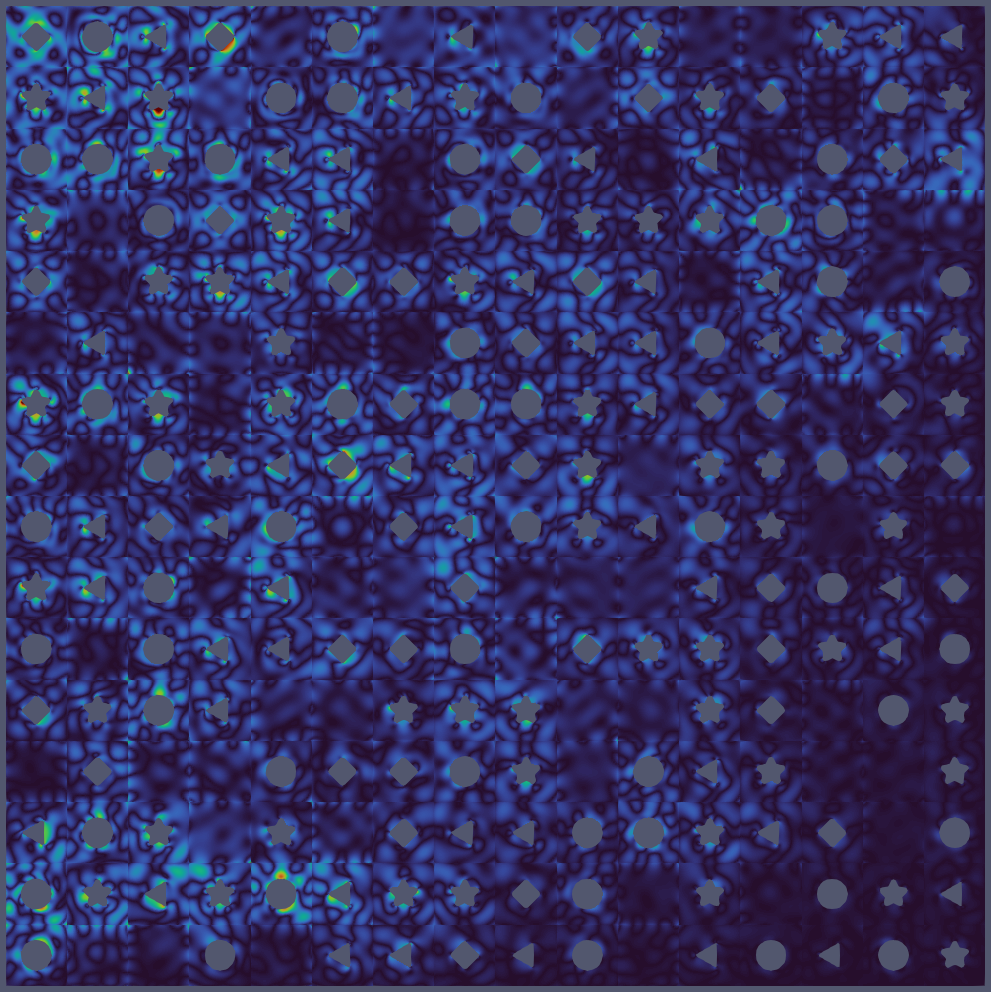}
            \addplot graphics[xmin=-0.1,xmax=16.1,ymin=-0.1,ymax=16.1]{\imagepath};

\end{groupplot}
            

\node[below = 1.5cm of my plots c1r1.south west,
    anchor=west,
] {(a) Velocity magnitude};
\node[below = 1.5cm of my plots c2r1.south west,
    anchor=west,
] {(b) Point-wise velocity error};
\node[below = 1.5cm of my plots c1r2.south west,
    anchor=west,
] {(c) Pressure};
\node[below = 1.5cm of my plots c2r2.south west,
    anchor=west,
] {(d) Point-wise pressure error};
\end{tikzpicture}
%
    \caption{Scaled-up ROM prediction of a $16\times16$ array at $\Re=25$ with basis sizes $R_{u,r}=R_{p,r}=Z_r=40$.}
    \label{fig:global-sol}
\end{figure}
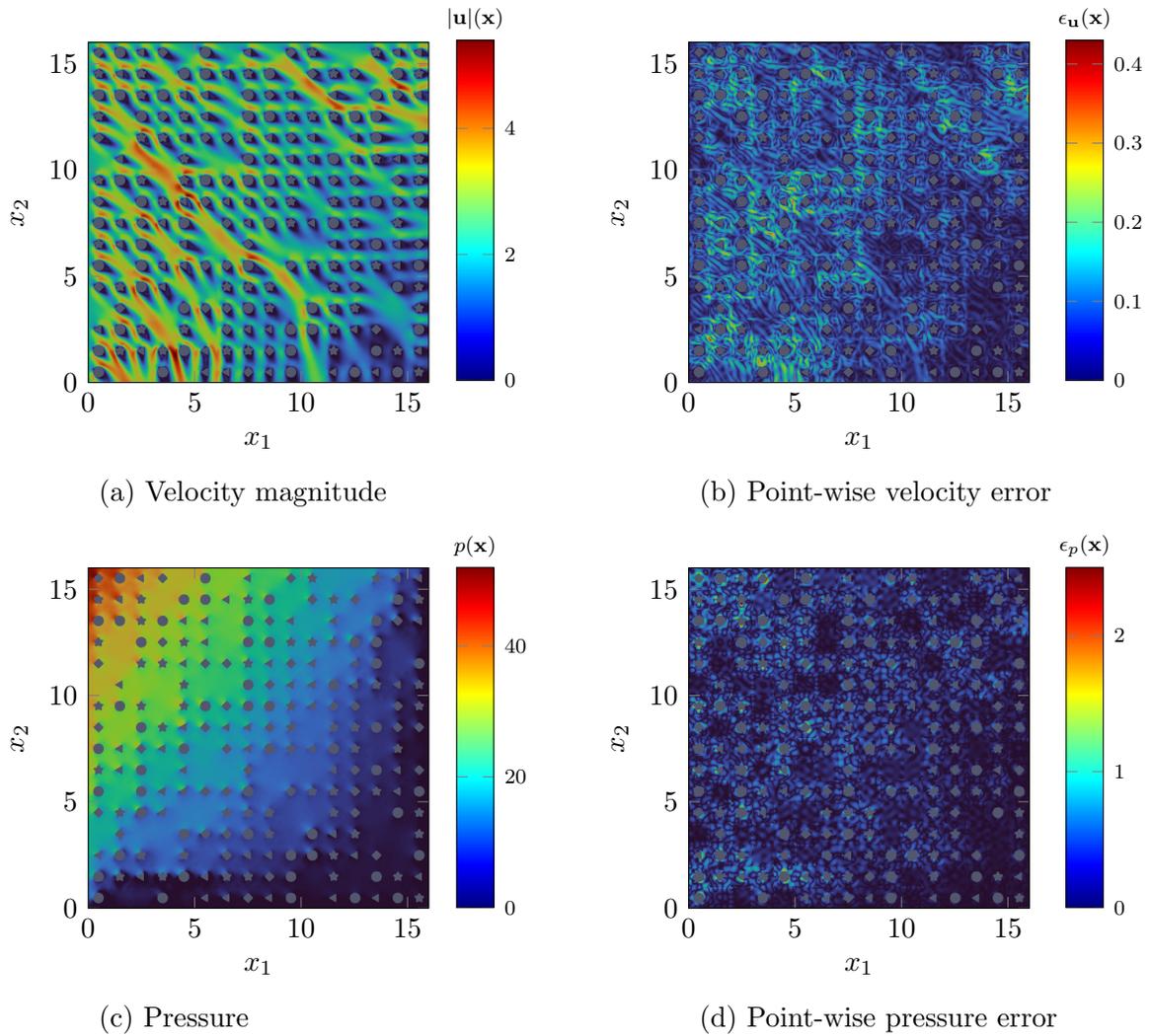
The global ROM (\ref{eq:rom-gov}) is built based on POD bases obtained in Section~\ref{subsec:bases}.
Figure~\ref{fig:global-sol} shows a scaled-up prediction of the global ROM for $16\times16$ array of objects,
with inflow velocity $\bg_{di}=(1.5, -0.8)$.
The dimensions of the bases are set to be $R_{u,r}=R_{p,r}=Z_r=40$ for all reference components, having $120$ basis vectors in total.
The tensorial approach (\ref{eq:rom-adv-tensor}) is used for nonlinear advection term,
though EQP (\ref{eq:eqp}) also produced almost identical results.
In Figure~\ref{fig:global-sol}~(a), large scale flow features are observed as empty components are combined together into large empty regions.
These features could not be observed in sample generation stage with $2\times2$ arrays.
Furthermore, flow tends to accumulate over these large empty regions,
increasing the velocity by a factor of about 5 compared to the training data.
Nonetheless, ROM was able to predict velocity with only $2.2\%$ relative error.
Figure~\ref{fig:global-sol}~(b) also shows point-wise absolute error of velocity.
Overall, the point-wise error is proportional to the velocity magnitude, up to a maximum of $10\%$ of its value.
\par
Similarly, ROM predictions for pressure also achieved $\sim0.9\%$ relative error.
In Figure~\ref{fig:global-sol}~(c), pressure decreases along the flow direction mainly across objects,
while remaining somewhat flat on empty regions.
In Figure~\ref{fig:global-sol}~(d), point-wise pressure error tends to scale with the magnitude of the pressure gradient,
having lower values in empty regions.
While there are localized high errors around the surface of objects,
the overall point-wise error remains less than $10\%$ of the pressure.
\par
\begin{figure}[tbph]
    \begin{tikzpicture}[font=\small,]
    \begin{groupplot}[
        group style={
            group name = my plots,
            group size= 2 by 2,
            xlabels at =edge bottom,
            horizontal sep=2cm,
            vertical sep=2.2cm,
        },
        name=chung,
    ]    
\pgfplotsset{set layers=standard}%

        \nextgroupplot[
            height = 0.45\textwidth,
            width = 0.5\textwidth,
            xlabel={Domain size $L^2$},
            ylabel={Assembly time ($s$)},
            tick scale binop ={\times},
            xmode=log, ymode=log,
            legend style={
                font=\small,
                draw=none, fill=none,
                at={(rel axis cs: 0., 1.0)},
                anchor=south west,
                nodes={scale=1.0},
                legend cell align={left},
                legend columns=4,
                /tikz/every even column/.append style={column sep=0.5cm},
            },
            legend image post style={mark options={solid, scale=1.0, fill=white, line width=1.0}},
        ]

            \addplot+ [
                line width=1.0,
                solid,
                mark=*,
                mark options={fill=white,},
                blue,
                error bars/.cd, y dir=both, y explicit,
            ]
            table [
                x expr=\thisrowno{0}^2, y index=1,
                y error minus expr=\thisrowno{1} - \thisrowno{4},
                y error plus expr=\thisrowno{5} - \thisrowno{1},
            ]{./data/scaling_fom_assemble_tensor.txt};

            \addplot+ [
                line width=1.0,
                solid,
                mark=*,
                mark options={fill=white,},
                red,
                error bars/.cd, y dir=both, y explicit,
            ]
            table [
                x expr=\thisrowno{0}^2, y index=1,
                y error minus expr=\thisrowno{1} - \thisrowno{4},
                y error plus expr=\thisrowno{5} - \thisrowno{1},
            ]{./data/scaling_rom_assemble_tensor.txt};

            \addplot+ [
                line width=1.0,
                dashed,
                mark=*,
                mark options={solid, fill=white,},
                red,
                error bars/.cd, y dir=both, y explicit,
            ]
            table [
                x expr=\thisrowno{0}^2, y index=1,
                y error minus expr=\thisrowno{1} - \thisrowno{4},
                y error plus expr=\thisrowno{5} - \thisrowno{1},
            ]{./data/scaling_rom_assemble_eqp.txt};

            \legend{FOM, {ROM, tensor}, {ROM, EQP}}

        \nextgroupplot[
            height = 0.45\textwidth,
            width = 0.5\textwidth,
            xlabel={Domain size $L^2$},
            ylabel={Solution time ($s$)},
            tick scale binop ={\times},
            xmode=log, ymode=log,
        ]
        
            \addplot+ [
                line width=1.0,
                solid,
                mark=*,
                mark options={fill=white,},
                blue,
                error bars/.cd, y dir=both, y explicit,
            ]
            table [
                x expr=\thisrowno{0}^2, y index=1,
                y error minus expr=\thisrowno{1} - \thisrowno{4},
                y error plus expr=\thisrowno{5} - \thisrowno{1},
            ]{./data/scaling_fom_solve_tensor.txt};

            \addplot+ [
                line width=1.0,
                solid,
                mark=*,
                mark options={fill=white,},
                red,
                error bars/.cd, y dir=both, y explicit,
            ]
            table [
                x expr=\thisrowno{0}^2, y index=1,
                y error minus expr=\thisrowno{1} - \thisrowno{4},
                y error plus expr=\thisrowno{5} - \thisrowno{1},
            ]{./data/scaling_rom_solve_tensor.txt};

            \addplot+ [
                line width=1.0,
                dashed,
                mark=*,
                mark options={fill=white, solid},
                red,
                error bars/.cd, y dir=both, y explicit,
            ]
            table [
                x expr=\thisrowno{0}^2, y index=1,
                y error minus expr=\thisrowno{1} - \thisrowno{4},
                y error plus expr=\thisrowno{5} - \thisrowno{1},
            ]{./data/scaling_rom_solve_eqp.txt};

            \logLogSlopeTriangle{0.6}{0.1}{0.32}{1.2}{}{north}
 
        \nextgroupplot[
            height = 0.45\textwidth,
            width = 0.5\textwidth,
            xlabel={Domain size $L^2$},
            ylabel={Relative error ($\%$)},
            tick scale binop ={\times},
            xmode=log,
            ymin=0., ymax=6,
            xshift=4cm,
            legend pos=north west,
            legend style={
                font=\small,
                draw=none, fill=none,
                anchor=north west,
                nodes={scale=1.0},
                legend cell align={left},
                legend columns=2,
            },
        ]

            \addplot+ [
                line width=1.0,
                solid,
                mark=*,
                mark options={fill=white,},
                brown,
                error bars/.cd, y dir=both, y explicit,
            ]
            table [
                x expr=\thisrowno{0}^2, y expr=\thisrowno{1} * 1e2,
                y error minus expr=\thisrowno{1} - \thisrowno{4},
                y error plus expr=\thisrowno{5} - \thisrowno{1},
            ]{./data/scaling_rel_error_vel_tensor.txt};

            \addplot+ [
                line width=1.0,
                solid,
                mark=*,
                mark options={fill=white,},
                purple,
                error bars/.cd, y dir=both, y explicit,
            ]
            table [
                x expr=\thisrowno{0}^2, y expr=\thisrowno{1} * 1e2,
                y error minus expr=\thisrowno{1} - \thisrowno{4},
                y error plus expr=\thisrowno{5} - \thisrowno{1},
            ]{./data/scaling_rel_error_pres_tensor.txt};

            \addplot+ [
                line width=1.0,
                dashed,
                mark=*,
                mark options={fill=white, solid},
                brown,
                error bars/.cd, y dir=both, y explicit,
            ]
            table [
                x expr=\thisrowno{0}^2, y expr=\thisrowno{1} * 1e2,
                y error minus expr=\thisrowno{1} - \thisrowno{4},
                y error plus expr=\thisrowno{5} - \thisrowno{1},
            ]{./data/scaling_rel_error_vel_eqp.txt};

            \addplot+ [
                line width=1.0,
                dashed,
                mark=*,
                mark options={fill=white, solid},
                purple,
                error bars/.cd, y dir=both, y explicit,
            ]
            table [
                x expr=\thisrowno{0}^2, y expr=\thisrowno{1} * 1e2,
                y error minus expr=\thisrowno{1} - \thisrowno{4},
                y error plus expr=\thisrowno{5} - \thisrowno{1},
            ]{./data/scaling_rel_error_pres_eqp.txt};

            \legend{{Tensor, $\bu$}, {Tensor, $\bp$}, {EQP, $\bu$}, {EQP, $\bp$}}

  \end{groupplot}
\end{tikzpicture}
%
    \caption{Prediction performance between FOM, ROM with tensorial approach (\ref{eq:rom-adv-tensor}) and EQP (\ref{eq:eqp})
    in (a) assembly time, (b) computation time and (c) accuracy. Error bars indicate $95\%$ confidence interval over 100 test cases.}
    \label{fig:scaling}
\end{figure}
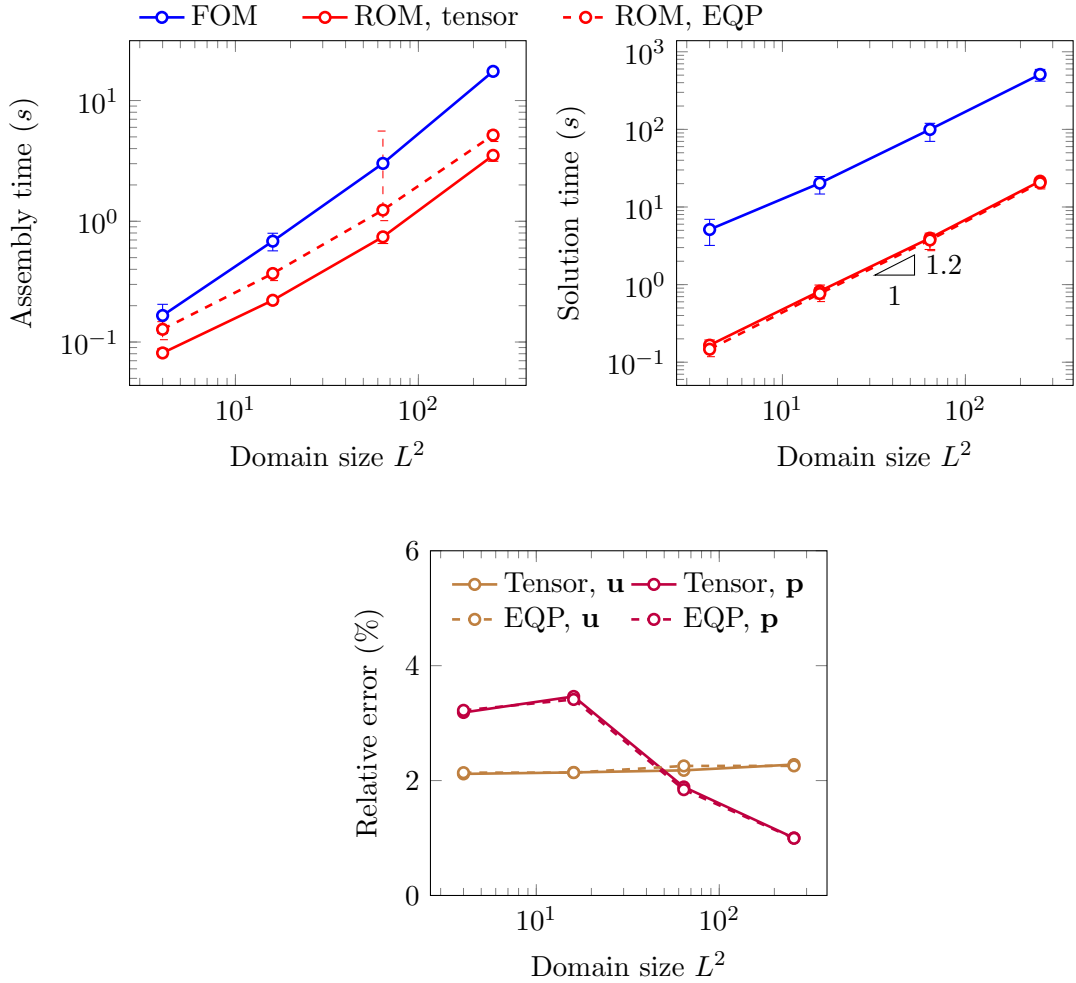
The global ROM is further demonstrated over different sizes of domain $\Omega=[0, L]^2$, from $L=2$ to $L=16$.
Same basis sizes $R_{u,r}=R_{p,r}=Z_r=40$ are used for all domain size.
At each domain size, similar to sample generation,
100 different random arrays are chosen with random inflow velocity,
in order to establish statistical significance.
Figure~\ref{fig:scaling} shows prediction performance of FOM and ROM over different sizes of domains.
Both FOM (\ref{eq:weak-fom-gov}) and ROM (\ref{eq:rom-gov}) are solved via the Newton-Raphson method.
Each iteration of Newton-Raphson method is solved using \texttt{MUMPS},
which is a direct LU-factorization linear solver~\cite{Amestoy2001,Amestoy2019}.
In Figure~\ref{fig:scaling}~(a), both the tensorial approach (\ref{eq:rom-adv-tensor}) and EQP (\ref{eq:eqp})
scale better with the domain size compared to FOM, achieving more speed-up in large domains.
For computation time in Figure~\ref{fig:scaling}~(b),
ROM achieved about $23.7\times$ speed-up over all scales compared to FOM,
with both tensorial approach and EQP.
Relative errors in Figure~\ref{fig:scaling}~(c), however,
remain less than $4\%$ over all scales.
Pressure relative error starts to decrease with domain size.
This is presumably because the pressure magnitude increases with domain size, as pressure accumulates along flow as shown in Figure~\ref{fig:global-sol}~(c).
Unlike pressure, velocity relative error remains about $2.28\%$ over all scales.
Overall, the results in Figure~\ref{fig:global-sol} shows that component ROM is capable of robust scaled-up prediction with much cheaper computational cost.
\par
In Figure~\ref{fig:global-sol}, the tensorial approach and EQP method did not differ much.
Their performance is further investigated subsequently in Section~\ref{subsec:tensor-vs-eqp}.

\subsection{Choice of basis and supremizer}
A numerical experiment is performed to demonstrate the role of the supremizer in Section~\ref{subsec:supreme}.
The global ROMs are built with $R_{u,r}=R_{p,r}=50$ POD basis vectors
and different numbers of supremizer basis vectors $Z_r=20,30,40,50$.
Their accuracy is compared based on predictions for a flow past $16\times16$ array of objects at $\Re=25$.
\par
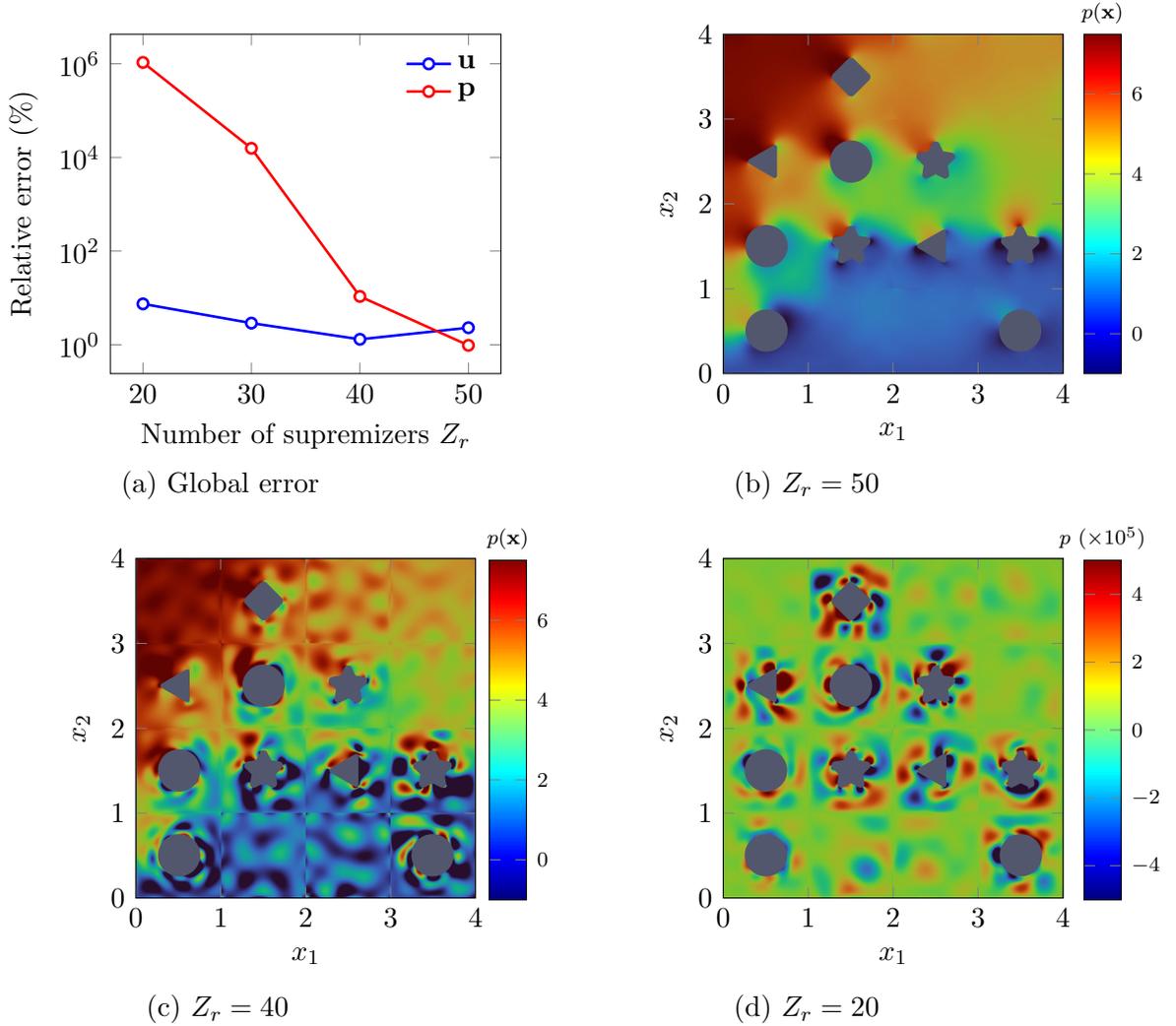
\begin{figure}[tbph]
    \begin{tikzpicture}[font=\small, spy using outlines={circle,black,magnification=6,size=1.5cm, connect spies}]
\pgfplotsset{set layers=standard}
    \begin{groupplot}[
        group style={
            group name = my plots,
            group size= 2 by 2,
            xlabels at =edge bottom,
            horizontal sep=3cm,
            vertical sep=2.5cm,
        },
        name=chung,
    ]    

        \nextgroupplot[
            height = 0.45\textwidth,
            width = 0.5\textwidth,
            xlabel={Number of supremizers $Z_r$},
            ylabel={Relative error ($\%$)},
            tick scale binop ={\times},
            ymode=log,
            legend pos=north east,
            legend style={
                font=\small,
                draw=none, fill=none,
                anchor=north east,
                nodes={scale=1.0},
                legend cell align={left},
            },
        ]

            \addplot+ [
                line width=1.0,
                solid,
                mark=*,
                mark options={fill=white,},
                blue,
            ]
            table [
                x index=0, y expr=\thisrowno{1} * 1e2,
            ]{./data/nsup.txt};

            \addplot+ [
                line width=1.0,
                solid,
                mark=*,
                mark options={fill=white,},
                red,
            ]
            table [
                x index=0, y expr=\thisrowno{2} * 1e2,
            ]{./data/nsup.txt};

            \legend{$\bu$, $\bp$}

        \nextgroupplot[
            height = 0.45\textwidth,
            width = 0.45\textwidth,
            ylabel={$x_2$},
            xlabel={$x_1$},
            tick scale binop ={\times},
            xmin = 0, xmax = 4,
            ymin = 0, ymax = 4,
            point meta min=-1, point meta max=7.5,
            colorbar, colormap/jet,
            colorbar style={
                font=\scriptsize,
                xticklabel pos=upper,
                scaled y ticks=false,
                /pgf/number format/precision=4,
                at={(rel axis cs: 2.8, 0.)}, anchor=south west,
                xlabel=$p(\bx)$,
            }
        ]
        
            \edef\imagepath{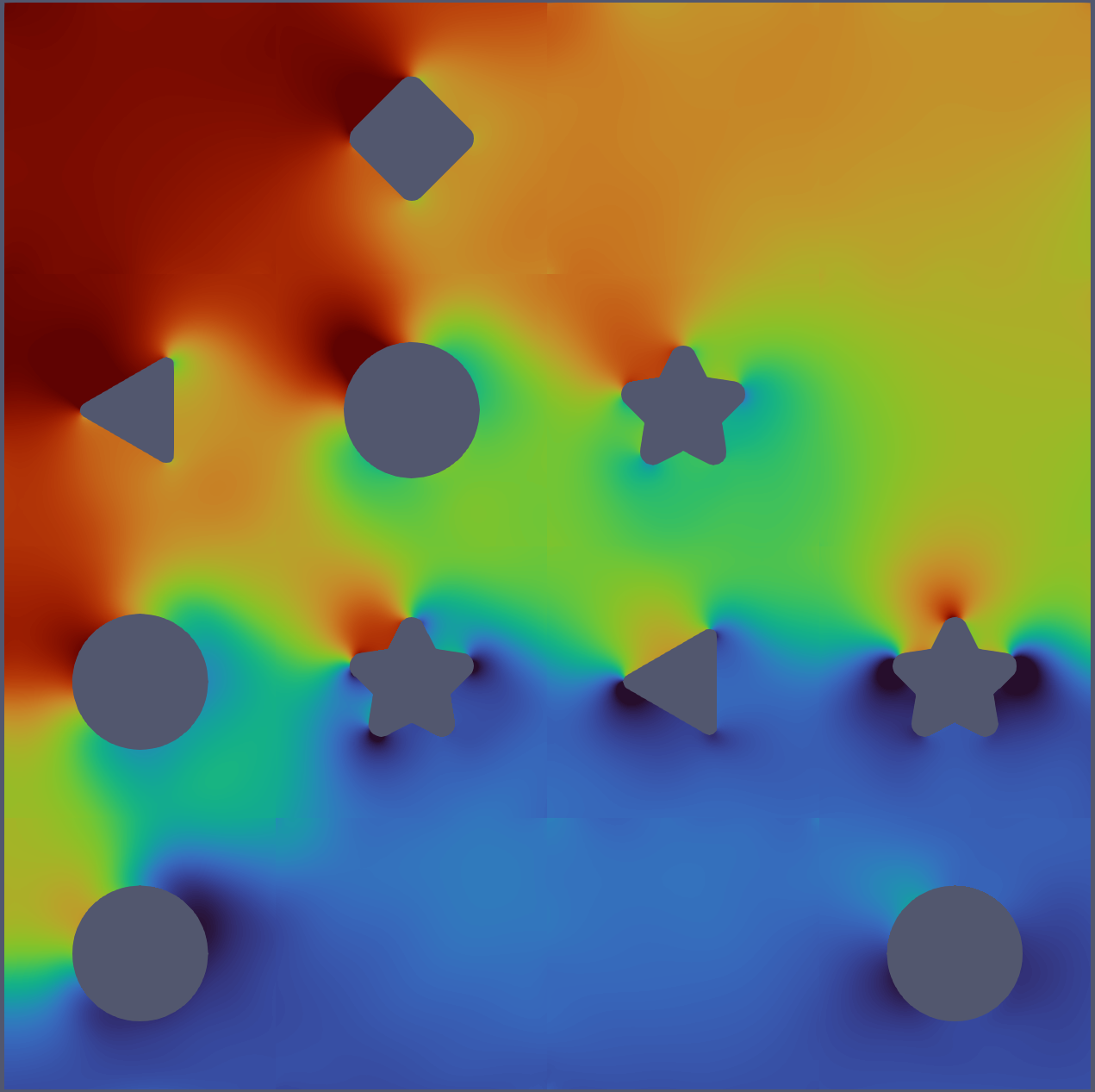}
            \addplot graphics[xmin=-0.,xmax=4.,ymin=-0.0,ymax=4.0]{\imagepath};

        \nextgroupplot[
            height = 0.45\textwidth,
            width = 0.45\textwidth,
            ylabel={$x_2$},
            xlabel={$x_1$},
            tick scale binop ={\times},
            xmin = 0, xmax = 4,
            ymin = 0, ymax = 4,
            point meta min=-1, point meta max=7.5,
            colorbar, colormap/jet,
            colorbar style={
                font=\scriptsize,
                xticklabel pos=upper,
                scaled y ticks=false,
                /pgf/number format/precision=4,
                at={(rel axis cs: 1.05, -1.55)}, anchor=south west,
                xlabel=$p(\bx)$,
            }
        ]
        
            \edef\imagepath{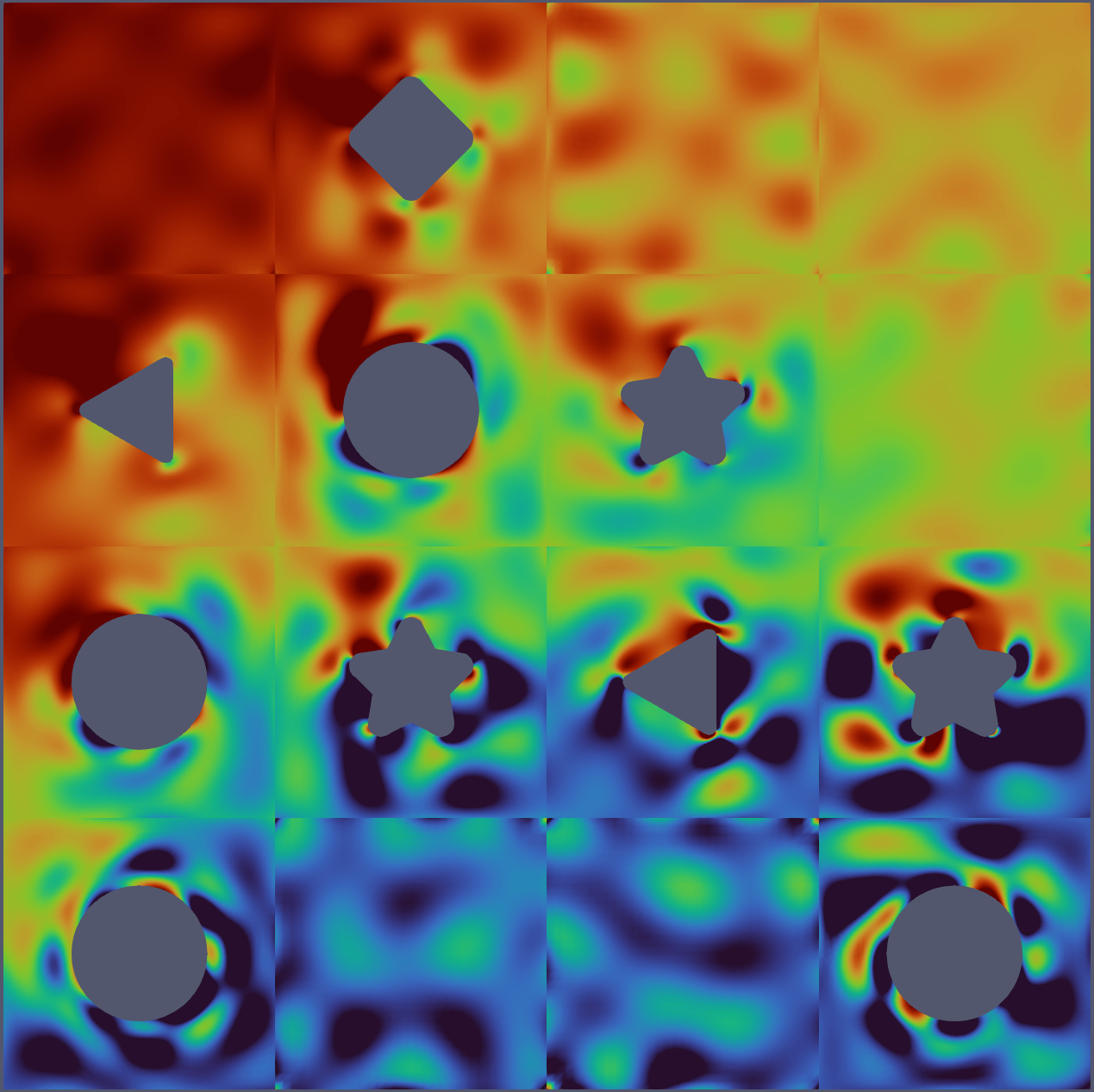}
            \addplot graphics[xmin=-0.,xmax=4.,ymin=-0.,ymax=4.]{\imagepath};

        \nextgroupplot[
            height = 0.45\textwidth,
            width = 0.45\textwidth,
            ylabel={$x_2$},
            xlabel={$x_1$},
            tick scale binop ={\times},
            xmin = 0, xmax = 4,
            ymin = 0, ymax = 4,
            point meta min=-5e5, point meta max=5e5,
            colorbar, colormap/jet,
            colorbar style={
                font=\scriptsize,
                xticklabel pos=upper,
                scaled y ticks=false,
                ytick={-4e5, -2e5, 0, 2e5, 4e5},
                yticklabels={{$-4$, $-2$, $0$, $2$, $4$}},
                /pgf/number format/precision=4,
                at={(rel axis cs: 2.8, -1.55)}, anchor=south west,
                xlabel=$p$ ($\times10^5$),
               x tick label style={
                   tick scale binop ={\times},
                   align=left,
                   /pgf/number format/fixed,
               },
            }
        ]
        
            \edef\imagepath{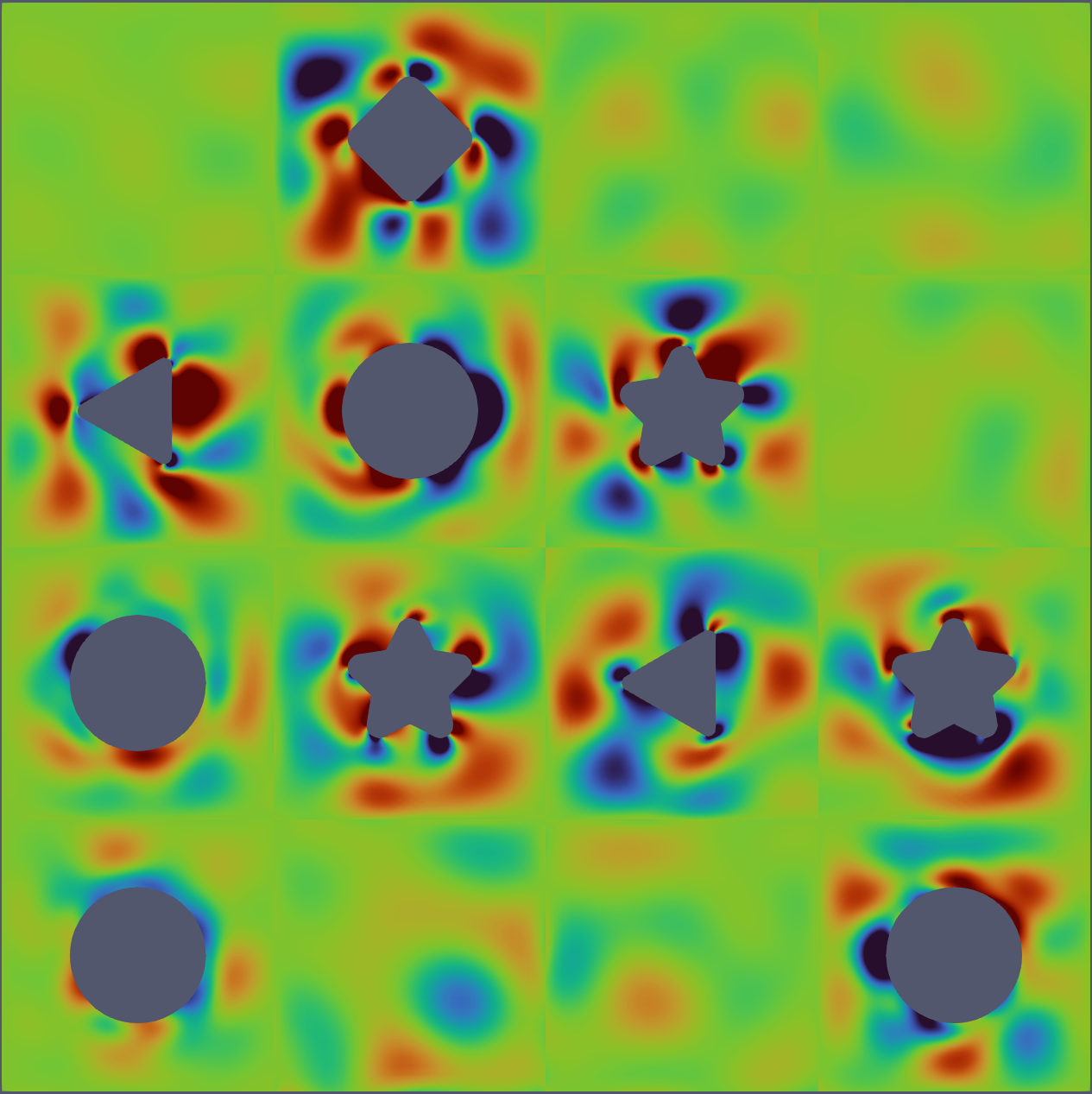}
            \addplot graphics[xmin=-0.,xmax=4.,ymin=-0.,ymax=4.]{\imagepath};

\end{groupplot}
            

\node[below = 1.5cm of my plots c1r1.south west,
    anchor=west,
] {(a) Global error};
\node[below = 1.5cm of my plots c2r1.south west,
    anchor=west,
] {(b) $Z_r=50$};
\node[below = 1.5cm of my plots c1r2.south west,
    anchor=west,
] {(c) $Z_r=40$};
\node[below = 1.5cm of my plots c2r2.south west,
    anchor=west,
] {(d) $Z_r=20$};
\end{tikzpicture}
%
    \caption{Effect of the supremizer on ROM prediction:
    (a) global error versus number of supremizers;
    magnified snapshot of ROM pressure prediction on $\bx\in[0, 4]^2$ with (b) $Z_r=50$;
    (c) $Z_r=40$;
    and (d) $Z_r=20$.
    For all unit components, $R_{u,r}=R_{p,r}=50$ POD basis vectors are used. \PR{Comment: Why the legends for color-bar are different in (b/c) and (d)? KC: due to its enormous magnitude, it didn't really visualize well in the same color range. I hoped to show the discontinuity over interfaces as well.}
    }
    \label{fig:supreme}
\end{figure}
Figure~\ref{fig:supreme}~(a) shows the global error of the ROM prediction according to the number of supremizers $Z_r$.
Even though supremizers increase the dimension for velocity subspaces represented by (\ref{eq:global-rom-test}),
increasing $Z_r$ only slightly improves the accuracy of the velocity from $7.5\%$ to $2\%$.
Increasing $Z_r$, however, significantly improves the accuracy of the pressure by orders of magnitude.
With $Z_r=R_{p,r}=50$, the ROM prediction recovers the pressure error as expected from the energy missing ratio (\ref{eq:pod-eps}) shown in Figure~\ref{fig:energy-ratio}~(b).
Though not reported here, simply increasing the pressure POD basis size $R_{p,r}$ rather aggravates the pressure error.
Figure~\ref{fig:supreme}~(b-d) shows pressure prediction with a decreasing number of supremizers.
With $Z_r=R_{p,r}=50$, the predicted pressure does not exhibit any spurious mode over interfaces or boundaries of objects.
Such spurious modes instantly appear as $Z_r$ decreases to 40.
With $Z_r=20$, the spurious modes diverge and dominate the entire domain.
\par
These results strongly indicate that the supremizers are necessary for stable pressure prediction in incompressible flow.
While velocity POD bases obtained from FOM snapshots may accurately represent velocity,
all basis vectors are divergence-free reflecting the incompressible nature from the snapshots.
The absence of compressible velocity components then leaves the pressure underdetermined,
thus introducing spurious pressure modes as observed.
The supremizer (\ref{eq:supreme}) provides an effective way to augment the velocity bases with compressible components,
thereby satisfying the necessary inf-sup condition~\cite{Babuvska1971,Brezzi1974,Ladyzhenskaya1963,Taylor1973}.
\par
We note that Chung \textit{et al.}~\cite{chung2024train} demonstrated a scaled-up prediction for incompressible Stokes flow
without supremizers. In their demonstration, velocity and pressure are treated as a combined variable $\bq = (\bu, \bp)$,
and the ROM representation (\ref{eq:linear-subspace}) is replaced with
\begin{equation}\label{eq:combined-pod}
    \bq_r = 
    \begin{pmatrix}
        \bu_r \\ \bp_r
    \end{pmatrix}
    \approx
    \bPhi_r\hbq_r
    =
    \begin{pmatrix}
        \bPhi_{u,r} \\ \bPhi_{p,r}
    \end{pmatrix}
    \hbq_r,
\end{equation}
where POD is performed over $\bq$-snapshots for $\bPhi_r$ as a whole, rather than $\bu$ and $\bp$ separately.
In this combined representation (\ref{eq:combined-pod}),
since the pressure is rather a linearly dependent variable of velocity,
the need for a supremizer is obviated.
However, this combined variable approach is not applicable to the steady Navier-Stokes equations,
as the nonlinear advection term necessitates the decoupling of velocity and pressure.
Though not reported here,
we observed in this study that extensions of the combined variable approach with the tensorial approach (\ref{eq:rom-adv-tensor}) or EQP (\ref{eq:eqp})
does not converge in Newton iteration.

\subsection{Comparison between the tensorial approach and EQP}\label{subsec:tensor-vs-eqp}

The tensorial approach (\ref{eq:rom-adv-tensor}) and EQP (\ref{eq:eqp}) are further compared with different basis sizes.
The comparison is conducted on 100 test cases of $16\times16$ arrays, each with random inflow velocity.
Bases sizes of velocity, pressure and supremizer are set uniform, i.e. $R_{u,r}=R_{p,r}=Z_r$.
The global ROM (\ref{eq:rom-gov}) is built using different basis sizes from $R_{u,r}=30$ to $R_{u,r}=80$, for both the tensorial approach and EQP.
\par
\begin{figure}[tbph]
    \begin{tikzpicture}[font=\small,]
    \begin{groupplot}[
        group style={
            group name = my plots,
            group size= 2 by 1,
            xlabels at =edge bottom,
            horizontal sep=2cm,
            vertical sep=1cm,
        },
        name=chung,
    ]        

    \nextgroupplot[
        height = 0.45\textwidth,
        width = 0.5\textwidth,
        xlabel={Basis dimension $R_{u,r}=R_{p,r}=Z_r$},
        ylabel={Solution time ($s$)},
        tick scale binop ={\times},
        xmode=log, ymode=log,
        xtick={30, 40, 50, 60, 70, 80},
        xticklabels={30, 40, 50, 60, 70, 80},
        legend style={
            font=\small,
            draw=none, fill=none,
            at={(rel axis cs: 0., 1.0)},
            anchor=south west,
            nodes={scale=1.0},
            legend cell align={left},
            legend columns=4,
            /tikz/every even column/.append style={column sep=0.5cm},
        },
        legend image post style={mark options={scale=1.0, fill=white, line width=1.0}},
    ]
    
        \addplot+ [
            line width=1.0,
            solid,
            mark=*,
            mark options={fill=white,},
            blue,
            error bars/.cd, y dir=both, y explicit,
        ]
        table [
            x index=0, y index=1,
            y error minus expr=\thisrowno{1} - \thisrowno{4},
            y error plus expr=\thisrowno{5} - \thisrowno{1},
        ]{./data/compare_rom_solve.tensor.txt};

        \addplot+ [
            line width=1.0,
            solid,
            mark=*,
            mark options={fill=white,},
            red,
            error bars/.cd, y dir=both, y explicit,
        ]
        table [
            x index=0, y index=1,
            y error minus expr=\thisrowno{1} - \thisrowno{4},
            y error plus expr=\thisrowno{5} - \thisrowno{1},
        ]{./data/compare_rom_solve.eqp.txt};

        \logLogReverseSlopeTriangle{0.65}{0.1}{0.75}{2.7}{}{south}
        \logLogSlopeTriangle{0.8}{0.1}{0.55}{2.2}{}{north}

        \legend{Tensor, EQP}

    \nextgroupplot[
        height = 0.45\textwidth,
        width = 0.5\textwidth,
        xlabel={Basis dimension $R_{u,r}=R_{p,r}=Z_r$},
        ylabel={Relative error ($\%$)},
        tick scale binop ={\times},
        xmode=log, ymode=log,
        ytick={0.3, 0.5, 1, 3, 5, 7},
        yticklabels={0.3, 0.5, 1, 3, 5, 7},
        xtick={30, 40, 50, 60, 70, 80},
        xticklabels={30, 40, 50, 60, 70, 80},
    ]

        \addplot+ [
            line width=1.0,
            solid,
            mark=*,
            mark options={fill=white,},
            blue,
            error bars/.cd, y dir=both, y explicit,
        ]
        table [
            x index=0, y expr=\thisrowno{1}*1e2,
            y error minus expr=(\thisrowno{1} - \thisrowno{4})*1e2,
            y error plus expr=(\thisrowno{5} - \thisrowno{1})*1e2,
        ]{./data/compare_rel_error_vel.tensor.txt};

        \addplot+ [
            line width=1.0,
            solid,
            mark=*,
            mark options={fill=white,},
            red,
            error bars/.cd, y dir=both, y explicit,
        ]
        table [
            x index=0, y expr=\thisrowno{1}*1e2,
            y error minus expr=(\thisrowno{1} - \thisrowno{4})*1e2,
            y error plus expr=(\thisrowno{5} - \thisrowno{1})*1e2,
        ]{./data/compare_rel_error_vel.eqp.txt};
        
        \logLogSlopeTriangle{0.8}{0.1}{0.4}{-3.5}{}{south}

  \end{groupplot}
\node[below = 1.5cm of my plots c1r1.south west,
    anchor=west,
] {(a)};
\node[below = 1.5cm of my plots c2r1.south west,
    anchor=west,
] {(b)};
\end{tikzpicture}
    \caption{Performance comparison between tensorial approach (\ref{eq:rom-adv-tensor}) and EQP (\ref{eq:eqp})
    in (a) computation time and (b) accuracy.}
    \label{fig:tensor-vs-eqp}
\end{figure}
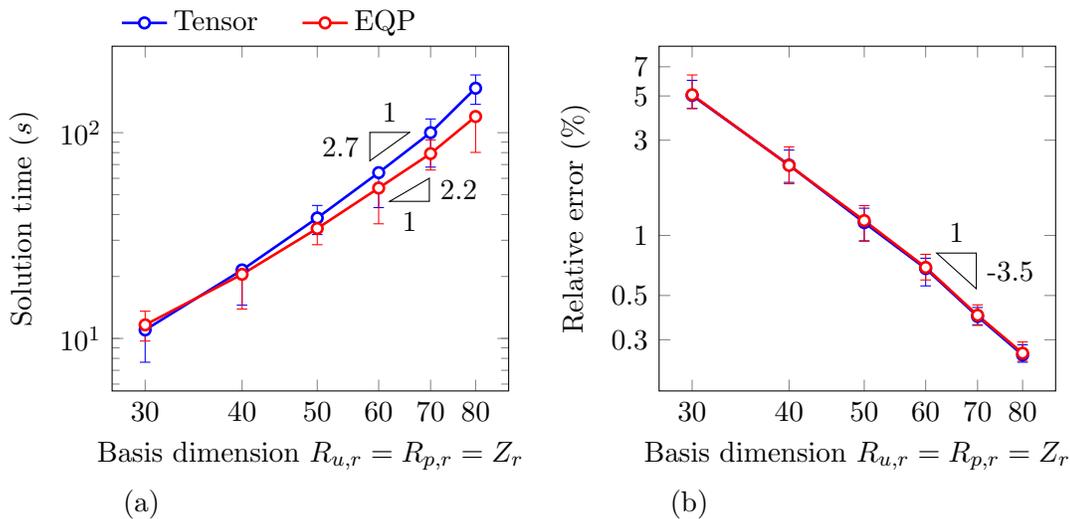
Figure~\ref{fig:tensor-vs-eqp} shows performance of the two nonlinear ROM approaches.
For basis size less then $50$, the performance of the two approaches are indistinguishable.
As the basis size increases, however,
computation time for EQP scales slightly better than tensorial approach, achieving more speed-up.
This did not come with compromise in accuracy, as both achieve same relative errors over all basis sizes.
For both approaches, relative error decays with basis size much faster than scaling of computation time,
showing the effectiveness of the nonlinear ROM approaches.
\par
While the comparison shows only a slightly better performance of EQP,
EQP has several advantages over tensorial approaches.
First of all, EQP is applicable to general nonlinear terms, unlike the tensorial approach which relies on the quadratic nature of the nonlinearity in the Navier-Stokes equations.
Also, computation time scaling of EQP is favorable for flow simulations at high Reynolds numbers,
where much larger basis sizes may be required to resolve the solution stably~\cite{tsai2023parametric}.
Another factor to be considered is the error threshold $\epsilon_{EQP}$ in (\ref{eq:eqp}),
which controls the number of empirical quadrature points and the corresponding accuracy.
While we followed the rule of thumb (\ref{eq:eqp-eps}),
the relative error in Figure~\ref{fig:tensor-vs-eqp}~(b) suggests that
the approximation error of empirical quadrature per (\ref{eq:eqp-eps}) is insignificant.
If larger error is permitted, EQP has another controlling factor for computation time by adjusting the error threshold,
while the tensorial approach does not.

\todo{PR:A brief discussion on how this method can be extended for higher Reynolds number cases should be added here. You can mention that the method has also been tested for other problems like backward step flow etc.\\
KC: it's in fact discussed at the end of the conclusion. added turbulent flows as another example there.}
\section{Conclusion}\label{sec:conclusion}

In this study, component reduced order modeling (CROM)~\cite{chung2024train} is extended to a nonlinear system, specifically, the steady Navier-Stokes equations.
Velocity and pressure snapshots are collected from flows past 2-by-2 arrays of five different unit component objects.
POD is performed for velocity and pressure of each reference component separately,
in order to obtain their respective linear subspace bases.
Velocity bases are augmented with compressible velocity components via a supremizer enrichment procedure~\cite{ballarin2015supremizer}.
A tensorial approach and EQP method~\cite{chapman2017accelerated} are introduced and compared
for tackling the ROM nonlinear advection operator.
The proposed method is demonstrated for flows past different sizes of arrays, up to 256 times larger than unit components.
For all scales, the global ROM achieved about $23$ times faster prediction with less than $4\%$ relative error.
A numerical experiment with the number of supremizers strongly indicates that
supremizer enrichment is essential for a stable pressure prediction.
The tensorial approach and EQP method are compared with different basis sizes,
and the advantages of EQP are discussed.
\par
The nonlinear ROM methods introduced in this study, in particular EQP method,
have a great potential for general nonlinear physics systems.
A natural extension is to couple with an advection-diffusion-reaction system for mass transfer problems. \TYL{For example, in CO$_2$ electrolysis, bicarbonate electrolytes are commonly used, where several different species are involved and reaction rates contain terms that are proportional to products of species concentrations \cite{moore2023simplified}.} \PR{In moisture sorption-diffusion by polymers, the diffusivity can be a function of relative humidity i.e. moisture concentration, and the water adsorption through molecular clustering can follow a power law~\citep{roy2022multi,castonguay2023modeling}.}
While the tensorial approach will be mainly used for advective terms,
EQP method can be employed for general nonlinear terms including Arrhenius-type reactions.
\par
More work remains to be done to uplift the applicability of the proposed method.
While the nonlinear ROM in this study does not involve interface or boundary terms,
such extension might be necessary for more complex and nonlinear physics problems.
This is particularly true for time-dependent hyperbolic conservation laws~\cite{Heath2012,Cockburn1999advdiff,Toro2013,Shahbazi2007,Ayuso2009,Houston2002},
where the solution interface condition is further constrained by nonlinear numerical fluxes.
Just as for linear interfacial operators, CROM can be employed without problem-specific interface handling
exploiting the fact that the interface condition is already handled at FOM level.
It is also well known that required basis sizes significantly increase for highly advection-dominated physics, e.g. \KC{high-Reynolds number turbulent flows or} a problem with Kolmogorov $n$-width decaying slowly.
For such cases, it would be of interest to investigate stabilization methods~\cite{tsai2023parametric},
parametrically local ROMs with multiple time-windows~\cite{copeland2022reduced,cheung2023local,Washabaugh2012,Amsallem2012}, and nonlinear manifold solution representation~\cite{kim2022fast, diaz2024fast, lauzon2024s}.

\section*{Acknowledgement}

This work was performed under the auspices of the U.S. Department of Energy
by Lawrence Livermore National Laboratory under contract DE-AC52-07NA27344
and was supported by Laboratory Directed Research and Development funding under project 22-SI-006.
LLNL-JRNL-870606.




\bibliographystyle{elsarticle-num} 
\bibliography{references}

%
%
%
%
\end{document}